\documentclass[10 pt, conference]{ieeeconf}
\pdfoutput=1

\IEEEoverridecommandlockouts                              % This command is only needed if
\usepackage{bbold}
\usepackage{color}
\usepackage{graphics} % for pdf, bitmapped graphics files
\usepackage{needspace}

\input{mysymbol.sty}
\usepackage{theorem}
\usepackage{cite}
\usepackage{amsmath}
\usepackage{mathtools}
\usepackage{multirow}
\usepackage{url}
\usepackage{graphics, times, amsfonts}
\usepackage{tikz, epic,eepic}
\usetikzlibrary{shapes,arrows}
\usepackage{pgfplots}
\usepackage{color}
\usepackage{epstopdf}
\usepackage{epsfig, amsbsy}
\usepackage{latexsym}
\usepackage{amscd, verbatim}
\usepackage{multirow}
\usepackage{booktabs}

\usepackage{tikz}
\usepackage{algorithm}
\usepackage{algorithmic}
\usepackage{subfiles}
\usepackage{bm}
\usepackage{mathtools}
\newtheorem{theorem}{Theorem}

\newtheorem{proposition}{Proposition}
\newtheorem{corollary}{Corollary}
\newtheorem{lemma}{Lemma}
\newtheorem{definition}{Definition}
\newtheorem{example}{Example}
{\itshape}{\rmfamily}

\newtheorem{remark}{Remark}

\newcommand{\Mat}[1]{\bm{\mathit{#1}}}

\def\forall{\text{for all\ }}

\title{\LARGE 
{Average Submodularity of Maximizing Anticoordination in Network Games}}
\author{Soham Das and Ceyhun Eksin % 
\thanks{Soham Das and Ceyhun Eksin are with the Industrial and Systems Engineering Department, Texas A\&M University, College Station, TX 77843. E-mail:{\tt\small  \; soham.das@tamu.edu; eksinc@tamu.edu.} This work was supported by NSF ECCS-1953694.%%
}}

\usepackage{graphicx}
\graphicspath{ {./images/} }
\usepackage{amssymb}

\begin{document}

\maketitle

\begin{abstract}
% In this paper we show that MAC is approximately submodular, hence a greedy approach to selecting control agents guarantees fractionally sub-optimal solutions for the line network. Moreover, we provide a special scenario for bipartite graphs that shows that the violation of submodularity can be arbitrarily bad, highlighting limitations of our provided guarantees in such worst case situations, thereby indicating that greedy approaches may not be always desirable in attaining anti-coordination. We further provide encouraging computational results which strongly suggest the effectiveness of greedy selection in practical scenarios.
%In this paper, we show that the maximum anti-coordination problem (MAC) is approximately submodular on line networks. Utilizing that result, we derive a performance guarantee for greedy agent selection algorithm for MAC. Finally, we provide a computational study to show the effectiveness of greedy node selection strategies to solve MAC on general bipartite networks.
% approximate submodularity can be arbitrarily bad, highlighting that greedy approaches will not be desirable in attaining anti-coordination.
We consider the control of decentralized learning dynamics for agents in an anti-coordination network game. In the anti-coordination network game, there is a preferred action in the absence of neighbors' actions, and the utility an agent receives from the preferred action decreases as more of its neighbors select the preferred action, potentially causing the agent to select a less desirable action. The decentralized dynamics that is based on the iterated elimination of dominated strategies converge for the considered game. Given a convergent action profile, we measure anti-coordination by the number of edges in the underlying graph that have at least one agent in either end of the edge not taking the preferred action. The maximum anti-coordination (MAC) problem seeks to find an optimal set of agents to control under a finite budget so that the overall network disconnect is maximized on game convergence as a result of the dynamics. We show that the MAC is submodular in expectation in dense bipartite networks for any  realization of the utility constants in the population. Utilizing this result, we obtain a performance guarantee for the greedy agent selection algorithm for MAC. Finally, we provide a computational study to show the effectiveness of greedy node selection strategies to solve MAC on general bipartite networks.
\end{abstract}
\section{Introduction}
% We consider a game of strategic substitutes played on a static network \cite{jackson2015games}. In a game with strategic substitutes, a player's incentive to play high value strategies decreases with an increase in the number of its neighbors playing a high strategy. We call this behaviour anti-coordination. 
Anti-coordination games can be used to study competition among firms \cite{bramoulle2007public,bramoulle2007anti}, public goods scenarios \cite{hirshleifer1983weakest}, free-rider behavior during epidemics \cite{bauch2013social,eksin2017disease}, and network security \cite{naghizadeh2016exit}. In each of these scenarios, there is a desired action for each agent, e.g., not taking the costly preemptive measures during a disease outbreak, not investing in insurance/protection etc., in the absence of other agents. When other agents are around, they can affect the benefits of the desired action, providing incentives for agents to switch. Here we consider networked interactions, where the actions of an agent are only affected by its neighbors (a subset of the population)---see  see \cite{parise2020analysis,galeotti2010network} for more details on network games. Despite the peer influences, some of the agents may continue to take the individually desired action, endangering their peers and the rest of the population. That is, the rational behavior can lead to the failure of anti-coordination in the population, when anti-coordination is desirable for the well-being and safety of the system as a whole. 

In such scenarios, we can envision the existence of a central coordinator with the goal to induce behavior that supports the well-being of the society. Here, we consider one such mechanism where the centralized coordinator intervenes by controlling a few agents in the network to incentivize anti-coordination among agents that repeatedly take actions to maximize individual payoffs. In particular, we consider decentralized learning dynamics inspired by the iterated elimination of dominated strategies \cite{menache2011network}. A dominated strategy is an action that cannot be preferred under any circumstance. Thus a dominated strategy cannot be a rational action. In the learning dynamics considered here, agents eliminate the dominated actions by evaluating their individual utilities under the worst and best possible action profiles of their neighbors (Section \ref{sec_learn}). The information about the elimination of an action by an agent can lead to a cascade of updates by other agents in the network. Indeed, we showed in \cite{eksin2020control} that such dynamics will converge in finite time and eliminate all dominated strategies of the anti-coordination network game considered here. Given such adaptive behavior of agents, the centralized player can steer the convergent action profile toward socially desirable outcomes by controlling the actions of a few players during the learning phase. 

In our setting, we define the goal of the central coordinator as to maximize anti-coordination between connected pairs of agents upon convergence of the behavior (Section \ref{sec_MAC}). The MAC problem is combinatorial, involving selection of a subset of the agents in the population. Indeed, we show that the MAC is NP-Hard to solve in general graphs (Theorem \ref{thm_hard}). 
This motivates us to consider a computationally tractable greedy selection protocol for solving MAC, where at each selection epoch the agent that yields the highest number of anti-coordinating edges at the convergent action profile is added to the control set until the control budget is reached (Section \ref{sec_MAC}). We show that the MAC problem is monotone and submodular in expectation almost surely in dense bipartite networks (Theorem \ref{thm_monotonic} and \ref{prop_avgsubmodular}). Together, these results imply that the worst case performance of the greedy selection protocol is bounded by a fraction of the optimal solution (Theorem \ref{thm_greedybound}). Numerical experiments show that the greedy selection protocol provides near-optimal results on average for general bipartite networks. 

This work is most closely related to the following intervention mechanisms in games that aim to improve efficiency: nudging \cite{xiao2011intervention,guers2013informational,riehl2018incentive}, influence maximization \cite{kempe2003maximizing}, and seeding in advertising \cite{balcan2012minimally}. All of these approaches aim to determine the emerging action profile resulting from an adaptive learning dynamics under repeated game play by either providing incentives or suggestions of ``good'' behavior to agents or by directly controlling a set of agents, as we do in this paper. Similar themes have been investigated in the control and optimization in networks literature \cite{bullo2009special}, flow of opinion dynamics in social networks \cite{bauso2016opinion,bolzern2019opinion}, and in the game theoretic control of multi-agent networked systems \cite{liu2019game}. Traditionally, a lot of work in this area is focused on consensus seeking or enhancing coordination among self interested agents with access to full or partial information. Here, we aim to maximize anti-coordination instead of maximizing social welfare or reach for a consensus. That is, we ask if selfish agents can successfully use local interaction to reach maximum polarization in strategies, and if not, can we identify a set of agents whose actions when controlled would lead to the same global objective? Other forms of intervention mechanisms involve financial incentives in the form of taxations or rewards \cite{brown2017studies}, and information design \cite{tavafoghi2017informational,sezer2021social}. These mechanisms do not consider repeated game play, and instead focus on improving the efficiency of Nash equilibria. Here our control selection policy is dependent on the adaptive behavior of agents. Lastly, we considered a similar MAC problem for the same anti-coordination network game in \cite{eksin2020control}. This work provides a performance guarantee for the greedy selection algorithm when agents interact over a dense bipartite network and there is a control budget.

% In Sections \ref{II} and \ref{III}, we first motivate and provide necessary definitions for the decentralized learning dynamics which govern the decision making of self interested agents who anti-coordinate in the network game. We motivate our setting by providing two examples. We talk about finite time convergence guarantees for the said dynamics and show that for special network structures, dynamics converge rapidly. 

% \red{We then begin our analysis in Section \ref{sec_MAC} formally defining the MAC optimization problem and by establishing that MAC is NP-Hard to solve in general graphs. This motivates us to look at approximate algorithms for solving MAC which may have provable performance guarantees. In particular, we explore the characteristics of a greedy agent selection algorithm for MAC.  We start Section V by formally defining the greedy approach, and defining submodularity for a set function. 

We obtain the performance guarantee for the greedy selection by showing the monotonicity and submodularity of the MAC problem in well-behaved instances given the learning dynamics. A well-behaved instance is one where there are no extremely insensitive agents, i.e. agents that are stubborn and play their chosen strategies no matter what their neighbors are playing. We capture this insensitivity to peer influence as having a very low learning constant. In such an instance, agents playing opposite strategies are concentrated on either side of the bipartite network. That is, if we control agents only on one side of the bipartite network, only agents on this side take the undesirable action as a result of the learning dynamics. 
% In the well-behaved scenario, agents playing opposite strategies are concentrated on either side of the bipartite network. 
This property makes anti-coordination akin to two parallel diffusion processes unfolding on both sides of the bipartite network. Each side reinforces the decision making of the other side, each side further polarizing the other side. We show that for dense networks, every instance is almost surely a well-behaved instance when we restrict the control set to a single side. 

Given the characteristics of the dynamics in well-behaved scenario, our approach to showing the submodularity of MAC follow ideas similar to the ones utilized in \cite{kempe2003maximizing,mossel2010submodularity}.
% Submodularity refers to a diminishing returns property.
% We endeavor to show in Theorem \ref{prop_avgsubmodular} that for dense bipartite networks, submodularity will hold in expectation. We introduce the notion of a well-behaved instance for MAC. In a well behaved instance, agents playing opposite strategies are concentrated on either side of the bipartite network. This property makes anti-coordination almost akin to two parallel diffusion processes unfolding on both sides of the bipartite graph. Each side reinforces the decision making of the other side, each side further polarizing the other side. We show that for dense networks, every instance is almost surely a well-behaved instance, and we shall exploit these behavioral characteristics of well behaved instances in our proofs. 
In \cite{kempe2003maximizing}, Kempe et al. talk about  an influence diffusion process (such as viral marketing),  where the process is initiated from a set of active nodes and spreads further. They show that if the dynamics are the linear threshold type, then the expected value of the total influence (the set of all active nodes) in the social network on process convergence is submodular. In \cite{mossel2010submodularity}, Mossel and Roch tackle the more general problem. They establish that any influence diffusion process, under certain assumptions on the neighbor influence function (monotonicity and submodularity) can be shown to be submodular. 

While our dynamics are  different, we are still able to utilize similar ideas to show submodularity (Theorem \ref{prop_avgsubmodular}). We begin by defining a function that measures one-step influence which encourages agents to adopt strategy level $0$ at each step of the dynamics. We establish that the one-step influence function is monotonic and submodular almost surely for dense networks (Proposition \ref{thm_submodularinfluence}). Next we describe the stochastic process of the set of agents taking the undesirable action ($0$) when the dynamics are seeded by a set of agents who are controlled. We show that the distribution of agents choosing to play strategy level $0$ on convergence of dynamics remains intact if instead of controlling the entire set from the get-go we break it into a set of smaller subsets (a partition) and control each set one by one in stages. This property allows us to provide an alternate equivalent description of the selection process based on a selection rule (Definition \ref{def_antisense} and Lemma \ref{lem_antisense}). %In Definition \ref{def_antisense} we formalize the selection rule and show that it is distribution preserving in Lemma \ref{lem_antisense}. 
The proof of the submodularity result entails designing such a {\it coupling} between the actual (greedy) selection method and another equivalent selection method for which we can show the diminishing returns property (submodularity). 
\section{Network anti-coordination game} \label{sec_network}

We consider a graph $\mathcal{G}(V,E)$ where the set of vertices $V=\{1,\dots,n\}$ represent the agents, and the set of edges $E$ represent interactions between the agents. Each agent takes an action $a_i \in [0,1]$ to maximize its utility function
\begin{equation} \label{eq_utility}
    u_i(a_i, a_{n(i)}) = a_i\big(1-c_i\sum_{j\in n(i)} a_j\big)
\end{equation}
where $n(i):=\{j:(i,j) \in E\}$ represents agent $i$'s neighbors, $a_{n(i)}:=\{a_j\}_{j\in n(i)}$ denotes the actions of agent $i$'s neighbors, and $0\le c_i<1$ is a constant. %The game above is a game of strategic substitutes. 
% Let us take $A_i$, the action space available to agent $i$, to be a complete lattice with an associated partial order $\geq_i$ for all $i\in V$. Then $A_i$ is a complete lattice if we define $a\geq a^{'}$ if and only if $a_i \geq a_i ^{'}$ for all $i$, and for any subset $S\subset A$, we can define the infimum and supremum naturally as $inf(S)=(inf_i\{a_i|a\in S\})_i$ and $sup(S)=(sup_i\{a_i|a\in S\})_i$. Then, a game is one of strategic substitutes if, for all $a_i \geq_i a_{i}'$ and $a_{-i}\geq_{-i}a_{-i}'$ 
% \begin{align}\label{eq_substitutes}
%     u_i(a_i,a_{-i})-u_i(a_{i}',a_{-i}) 
%     \leq u_i(a_i,a_{-i}')-u_i(a_{i}',a_{-i}') 
% \end{align}

\begin{figure*}
\centering
% !TEX root = ../Control in state-based anti-coordination network games.tex

% \usetikzlibrary{matrix,arrows,decorations.pathmorphing}
% \usetikzlibrary{arrows,automata}

\usetikzlibrary{matrix}
\usetikzlibrary{arrows}
\usetikzlibrary{decorations.pathmorphing}
% \small{\textit{Situation 1}}
\tikzstyle{sets} = [ellipse, draw=black, inner sep=0pt, minimum size=0.5cm]
\tikzstyle{terminal} = [circle, draw=black, inner sep=0pt, minimum size=0.5cm]
\tikzstyle{dashcircle}=[circle,draw=red,dashed,inner sep=0pt, minimum size=0.7cm]
\tikzstyle{theta} = [circle, draw=white, inner sep=0pt, minimum size=0.5cm]
\tikzstyle{arrow} = [stealth-stealth, thick]
\tikzstyle{edge}=[ultra thin, draw=black!30]

\begin{tikzpicture}[x=0.4cm, y=0.3cm]
%Draw grid as reference mark
% \draw[step=1.5, color=black!10] (0,-5) grid (9,5);
%Nodes

\filldraw[] (2.5,3) node (terminal 5) [terminal,label=90:{}]{$\epsilon$} 
% ++ (0:0) node (dashcircle )[dashcircle]{}
++ (270:5) node (terminal 1) [terminal,color=black!100,ultra thick,label=270:{}]{$\epsilon$}
++ (60:5.773) node (terminal 6) [terminal,label=90:{}]{$\epsilon$}
++ (270:5) node (terminal 2)[terminal,color=black!100,ultra thick,label=270:{} ]{$\epsilon$}
% ++ (0:0) node (dashcircle) [dashcircle]{$\epsilon$}
++ (60:5.773) node (terminal 7)[terminal,label=90:{}]{$\epsilon$}
++ (270:5) node (terminal 3) [terminal,color=black!100,ultra thick,label=270:{}]{$\epsilon$}
++ (60:5.773) node (terminal 8) [terminal,label=90:{}]{$\epsilon$}
++ (270:5) node (terminal 4) [terminal,color=black!100,ultra thick,label=270:{}]{$\epsilon$};

% \draw[edge] (terminal 1) to (terminal 2) to (terminal 3) to (terminal 4) to (terminal 5) to (terminal 6) to (terminal 7) to (terminal 8) to (terminal 5) to (terminal 2) to (terminal 7) to (terminal 4) to (terminal 1) to (terminal 8) to (terminal 3) to (terminal 6) to (terminal 1);

\draw[edge] (terminal 1) to (terminal 5);
\draw[edge] (terminal 1) to (terminal 6);
\draw[edge] (terminal 1) to (terminal 7);
\draw[edge] (terminal 2) to (terminal 5);
\draw[edge] (terminal 2) to (terminal 8);
\draw[edge] (terminal 3) to (terminal 6);
\draw[edge] (terminal 3) to (terminal 7);
\draw[edge] (terminal 3) to (terminal 8);
\draw[edge] (terminal 4) to (terminal 5);
\draw[edge] (terminal 4) to (terminal 7);
\draw[edge] (terminal 4) to (terminal 8);

% Draw S0 and S1 and the rectangles
\draw[draw=red,dashed] (1.5,-4) rectangle ++(10.8,4);
% \node[draw, align =left, font=\footnotesize] at (15,-2){$S_0$};
\draw[draw=blue,dashed] (1.5,1) rectangle ++(10.8,4);
% \node[draw, align =left, font=\footnotesize] at (15,3){$S_1$};

\filldraw[] (16.5,3) node (terminal 13) [terminal]{$\epsilon$} 
% ++ (0:0) node (dashcircle )[dashcircle]{}
++ (270:5) node (terminal 9) [terminal,color=black!100,ultra thick]{$\epsilon$}
++ (60:5.773) node (terminal 14) [terminal]{$\epsilon$}
++ (270:5) node (terminal 10)[terminal,color=black!100,ultra thick,label=270:{} ]{$\epsilon$}
% ++ (0:0) node (dashcircle) [dashcircle]{}
++ (60:5.773) node (terminal 15)[terminal]{$\epsilon$}
++ (270:5) node (terminal 11) [terminal,color=red!100,ultra thick,label=270:{}]{$0$}
++ (60:5.773) node (terminal 16) [terminal]{$\epsilon$}
++ (270:5) node (terminal 12) [terminal,color=red!100,ultra thick]{$0$};

\draw[edge] (terminal 9) to (terminal 13);
\draw[edge] (terminal 9) to (terminal 14);
\draw[edge] (terminal 9) to (terminal 15);
\draw[edge] (terminal 10) to (terminal 13);
\draw[edge] (terminal 10) to (terminal 16);
\draw[edge] (terminal 11) to (terminal 14);
\draw[edge] (terminal 11) to (terminal 15);
\draw[edge] (terminal 11) to (terminal 16);
\draw[edge] (terminal 12) to (terminal 13);
\draw[edge] (terminal 12) to (terminal 15);
\draw[edge] (terminal 12) to (terminal 16);
\draw[draw=red,dashed] (15.5,-4) rectangle ++(10.8,4);
% \node[draw, align =left, font=\footnotesize] at (15,-2){$S_0$};
\draw[draw=blue,dashed] (15.5,1) rectangle ++(10.8,4);

% \filldraw[] (31,3) node (terminal 21) [terminal]{$\epsilon$} 
% ++ (270:5) node (terminal 17) [terminal,color=black!100,ultra thick]{$\epsilon$}
% ++ (60:5.773) node (terminal 22) [terminal]{$1$}
% ++ (270:5) node (terminal 18)[terminal,color=black!100,ultra thick,label=270:{} ]{$\epsilon$}
% % ++ (0:0) node (dashcircle) [dashcircle]{}
% ++ (60:5.773) node (terminal 23)[terminal]{$1$}
% ++ (270:5) node (terminal 19) [terminal,color=red!100,ultra thick,label=270:{}]{$0$}
% ++ (60:5.773) node (terminal 24) [terminal]{$1$}
% ++ (270:5) node (terminal 20) [terminal,color=red!100,ultra thick]{$0$};

% \draw[edge] (terminal 17) to (terminal 21);
% \draw[edge] (terminal 17) to (terminal 22);
% \draw[edge] (terminal 17) to (terminal 23);
% \draw[edge] (terminal 18) to (terminal 21);
% \draw[edge] (terminal 18) to (terminal 24);
% \draw[edge] (terminal 19) to (terminal 22);
% \draw[edge] (terminal 19) to (terminal 23);
% \draw[edge] (terminal 19) to (terminal 24);
% \draw[edge] (terminal 20) to (terminal 21);
% \draw[edge] (terminal 20) to (terminal 23);
% \draw[edge] (terminal 20) to (terminal 24);

% \draw[draw=red,dashed] (30,-4) rectangle ++(10.8,4);
% % \node[draw, align =left, font=\footnotesize] at (15,-2){$S_0$};
% \draw[draw=blue,dashed] (30,1) rectangle ++(10.8,4);
% \node[draw, align =left, font=\footnotesize] at (14.5,-6){$c_1=0.41,c_2=0.55,c_3=0.57,c_4=0.86,c_5=0.92,c_6=0.60,c_7=0.34,c_8=0.39$};
% \node[draw, align=left, font=\footnotesize] at (2.5,4.5){$5$};

% \node[draw, align =left, font=\footnotesize] at (6.5,7){$c_5=0.92,c_6=0.60,c_7=0.34,c_8=0.39$};

\node[draw, align =left, font=\footnotesize] at (5.5,9){$t=1$};
\node[draw, align =left, font=\footnotesize] at (19.5,9){$t=2$};
% \node[draw, align =left, font=\footnotesize] at (33.5,9){$t=3$};

\filldraw[] (2.5,6) node [font=\footnotesize]{$c_5$} 
% ++ (0:0) node (dashcircle )[dashcircle]{}
++ (0:3) node  [font=\footnotesize]{$c_6$}
++ (0:3) node  [font=\footnotesize]{$c_7$}
++ (0:3) node  [font=\footnotesize]{$c_8$};

\filldraw[] (2.5,-5.5) node [font=\footnotesize]{$c_1$} 
% ++ (0:0) node (dashcircle )[dashcircle]{}
++ (0:3) node  [font=\footnotesize]{$c_2$}
++ (0:3) node  [font=\footnotesize]{$c_3$}
++ (0:3) node  [font=\footnotesize]{$c_4$};
\end{tikzpicture}
% !TEX root = ../Control in state-based anti-coordination network games.tex

% \usetikzlibrary{matrix,arrows,decorations.pathmorphing}
% \usetikzlibrary{arrows,automata}

\usetikzlibrary{matrix}
\usetikzlibrary{arrows}
\usetikzlibrary{decorations.pathmorphing}
% \small{\textit{Situation 1}}
\tikzstyle{sets} = [ellipse, draw=black, inner sep=0pt, minimum size=0.5cm]
\tikzstyle{terminal} = [circle, draw=black, inner sep=0pt, minimum size=0.5cm]
\tikzstyle{dashcircle}=[circle,draw=red,dashed,inner sep=0pt, minimum size=0.7cm]
\tikzstyle{theta} = [circle, draw=white, inner sep=0pt, minimum size=0.5cm]
\tikzstyle{arrow} = [stealth-stealth, thick]
\tikzstyle{edge}=[ultra thin, draw=black!30]

\begin{tikzpicture}[x=0.4cm, y=0.3cm]
%Draw grid as reference mark
% \draw[step=1.5, color=black!10] (0,-5) grid (9,5);
%Nodes

\filldraw[] (2.5,3) node (terminal 5) [terminal,label=90:{}]{$\epsilon$} 
% ++ (0:0) node (dashcircle )[dashcircle]{}
++ (270:5) node (terminal 1) [terminal,color=black!100,ultra thick,label=270:{}]{$\epsilon$}
++ (60:5.773) node (terminal 6) [terminal,label=90:{}]{$1$}
++ (270:5) node (terminal 2)[terminal,color=black!100,ultra thick,label=270:{} ]{$\epsilon$}
% ++ (0:0) node (dashcircle) [dashcircle]{$\epsilon$}
++ (60:5.773) node (terminal 7)[terminal,label=90:{}]{$1$}
++ (270:5) node (terminal 3) [terminal,color=red!100,ultra thick,label=270:{}]{$0$}
++ (60:5.773) node (terminal 8) [terminal,label=90:{}]{$1$}
++ (270:5) node (terminal 4) [terminal,color=red!100,ultra thick,label=270:{}]{$0$};

% \draw[edge] (terminal 1) to (terminal 2) to (terminal 3) to (terminal 4) to (terminal 5) to (terminal 6) to (terminal 7) to (terminal 8) to (terminal 5) to (terminal 2) to (terminal 7) to (terminal 4) to (terminal 1) to (terminal 8) to (terminal 3) to (terminal 6) to (terminal 1);

\draw[edge] (terminal 1) to (terminal 5);
\draw[edge] (terminal 1) to (terminal 6);
\draw[edge] (terminal 1) to (terminal 7);
\draw[edge] (terminal 2) to (terminal 5);
\draw[edge] (terminal 2) to (terminal 8);
\draw[edge] (terminal 3) to (terminal 6);
\draw[edge] (terminal 3) to (terminal 7);
\draw[edge] (terminal 3) to (terminal 8);
\draw[edge] (terminal 4) to (terminal 5);
\draw[edge] (terminal 4) to (terminal 7);
\draw[edge] (terminal 4) to (terminal 8);

% Draw S0 and S1 and the rectangles
\draw[draw=red,dashed] (1.5,-4) rectangle ++(10.8,4);
% \node[draw, align =left, font=\footnotesize] at (15,-2){$S_0$};
\draw[draw=blue,dashed] (1.5,1) rectangle ++(10.8,4);
% \node[draw, align =left, font=\footnotesize] at (15,3){$S_1$};

\filldraw[] (16.5,3) node (terminal 13) [terminal]{$\epsilon$} 
% ++ (0:0) node (dashcircle )[dashcircle]{}
++ (270:5) node (terminal 9) [terminal,color=black!100,ultra thick]{$\epsilon$}
++ (60:5.773) node (terminal 14) [terminal]{$1$}
++ (270:5) node (terminal 10)[terminal,color=black!100,ultra thick,label=270:{} ]{$\epsilon$}
% ++ (0:0) node (dashcircle) [dashcircle]{}
++ (60:5.773) node (terminal 15)[terminal]{$1$}
++ (270:5) node (terminal 11) [terminal,color=red!100,ultra thick,label=270:{}]{$0$}
++ (60:5.773) node (terminal 16) [terminal]{$1$}
++ (270:5) node (terminal 12) [terminal,color=red!100,ultra thick]{$0$};

\draw[edge] (terminal 9) to (terminal 13);
\draw[edge] (terminal 9) to (terminal 14);
\draw[edge] (terminal 9) to (terminal 15);
\draw[edge] (terminal 10) to (terminal 13);
\draw[edge] (terminal 10) to (terminal 16);
\draw[edge] (terminal 11) to (terminal 14);
\draw[edge] (terminal 11) to (terminal 15);
\draw[edge] (terminal 11) to (terminal 16);
\draw[edge] (terminal 12) to (terminal 13);
\draw[edge] (terminal 12) to (terminal 15);
\draw[edge] (terminal 12) to (terminal 16);
\draw[draw=red,dashed] (15.5,-4) rectangle ++(10.8,4);
% \node[draw, align =left, font=\footnotesize] at (15,-2){$S_0$};
\draw[draw=blue,dashed] (15.5,1) rectangle ++(10.8,4);

% \filldraw[] (31,3) node (terminal 21) [terminal]{$\epsilon$} 
% ++ (270:5) node (terminal 17) [terminal,color=black!100,ultra thick]{$\epsilon$}
% ++ (60:5.773) node (terminal 22) [terminal]{$1$}
% ++ (270:5) node (terminal 18)[terminal,color=black!100,ultra thick,label=270:{} ]{$\epsilon$}
% % ++ (0:0) node (dashcircle) [dashcircle]{}
% ++ (60:5.773) node (terminal 23)[terminal]{$1$}
% ++ (270:5) node (terminal 19) [terminal,color=red!100,ultra thick,label=270:{}]{$0$}
% ++ (60:5.773) node (terminal 24) [terminal]{$1$}
% ++ (270:5) node (terminal 20) [terminal,color=red!100,ultra thick]{$0$};

% \draw[edge] (terminal 17) to (terminal 21);
% \draw[edge] (terminal 17) to (terminal 22);
% \draw[edge] (terminal 17) to (terminal 23);
% \draw[edge] (terminal 18) to (terminal 21);
% \draw[edge] (terminal 18) to (terminal 24);
% \draw[edge] (terminal 19) to (terminal 22);
% \draw[edge] (terminal 19) to (terminal 23);
% \draw[edge] (terminal 19) to (terminal 24);
% \draw[edge] (terminal 20) to (terminal 21);
% \draw[edge] (terminal 20) to (terminal 23);
% \draw[edge] (terminal 20) to (terminal 24);

% \draw[draw=red,dashed] (30,-4) rectangle ++(10.8,4);
% % \node[draw, align =left, font=\footnotesize] at (15,-2){$S_0$};
% \draw[draw=blue,dashed] (30,1) rectangle ++(10.8,4);
% \node[draw, align =left, font=\footnotesize] at (14.5,-6){$c_1=0.41,c_2=0.55,c_3=0.57,c_4=0.86,c_5=0.92,c_6=0.60,c_7=0.34,c_8=0.39$};
% \node[draw, align=left, font=\footnotesize] at (2.5,4.5){$5$};

% \node[draw, align =left, font=\footnotesize] at (6.5,7){$c_5=0.92,c_6=0.60,c_7=0.34,c_8=0.39$};

\node[draw, align =left, font=\footnotesize] at (5.5,7){$t=3$};
\node[draw, align =left, font=\footnotesize] at (19.5,7){$t=4$};

% \filldraw[] (2.5,6) node [font=\footnotesize]{$c_5$} 
% % ++ (0:0) node (dashcircle )[dashcircle]{}
% ++ (0:3) node  [font=\footnotesize]{$c_6$}
% ++ (0:3) node  [font=\footnotesize]{$c_7$}
% ++ (0:3) node  [font=\footnotesize]{$c_8$};

% \filldraw[] (2.5,-5.5) node [font=\footnotesize]{$c_1$} 
% % ++ (0:0) node (dashcircle )[dashcircle]{}
% ++ (0:3) node  [font=\footnotesize]{$c_2$}
% ++ (0:3) node  [font=\footnotesize]{$c_3$}
% ++ (0:3) node  [font=\footnotesize]{$c_4$};
\end{tikzpicture}
\caption{{We have a $4\times 4$ bipartite well-behaved instance for MAC, with learning constants for agents given as $c_1=0.41,c_2=0.55,c_3=0.57,c_4=0.86,c_5=0.92,c_6=0.60,c_7=0.34,c_8=0.39$. See that all $c_i, i\in V$ satisfy $c_i \geq 1/d_i$ where $d_i$ is the node degree of agent $i$ in this example. Sets $S_0$ and $S_1$ are marked in red and blue boxes. At time $t=1$, all agents are undecided. We initiate dynamics by applying control to agents in set $S_0$ (shown in red) at $t=2$. The dynamics in \eqref{eq_updates} converge at $t=3$. No agent, following the dynamics would change their strategy beyond $t=3$.}}
\label{fig:bipartite1}
\end{figure*}
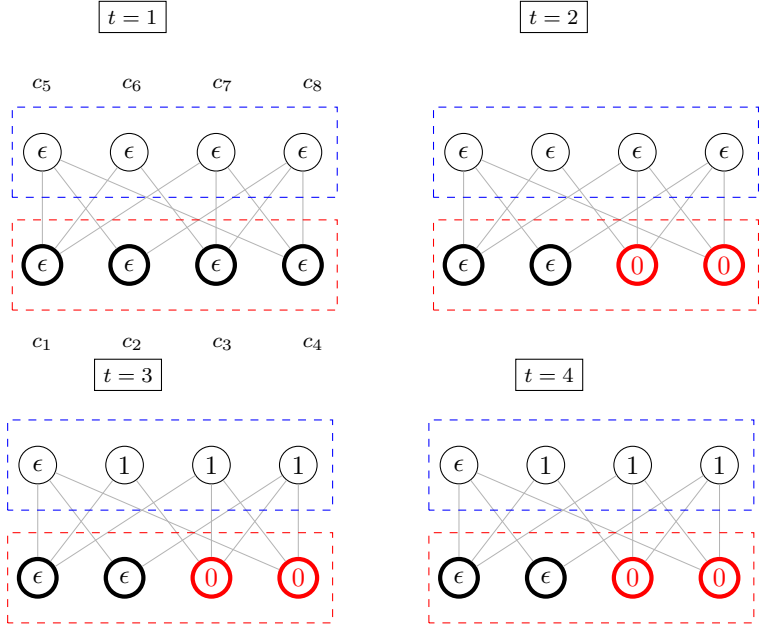

% The utility function \eqref{eq_utility} satisfies \eqref{eq_substitutes}. 
The utility function above captures scenarios where agent $i$ has a preferred action (action 1) but its incentive to choose this action decreases as more of its neighbors choose the preferred action. The decrease in incentive per neighbor is proportional to the constant $c_i$. The constant represents the sensitivity of agent $i$'s utility to its neighbor's actions.

The network game is represented by the tuple $\Gamma = \{V, A,\{u_i\}_{i\in V}\}$ where $A=[0,1]^{n}$ is the set of actions available to all players. 
% Note here that while the whole $[0,1]$ space is available to the agent to choose its action from, utility maximization enables either of the two possibilities: select $1$ when the quantity inside parenthesis in \eqref{eq_utility} is non-negative, or select $0$ otherwise. 

Below we provide scenarios that can be captured by the network game $\Gamma$ with payoffs given by \eqref{eq_utility}. 

\begin{example}[Disease spread on networks \cite{eksin2017disease}] \label{example1}
Consider a bipartite graph $\mathcal{G}_B$ where agents on one side are sick and the other side healthy \cite{eksin2017disease}. The edges in $\mathcal{G}_B$ represent the network of interaction between them. The disease spreads when agents on either end of a interaction link do not follow healthcare protocols, such as wearing masks, vaccinating, social distancing etc. We model this using action level $1$ for the agent (the easier/preferred action). Action level $0$ represents following epidemic mitigation protocols (the costlier action), and all actions between $0$ and $1$ represents the relative importance given to disease prevention measures. When we have an agent playing $0$ on the end of an interaction link, we have deactivated a disease transmission pathway in society. 
%\red{Given an interaction graph, the players are exposed to the strategies of their neighbors in the previous time step, and players simultaneously best-respond to that information.}\blue{CE: At this point the reader does not know anything about the dynamics.} 
The utilities of agents in the epidemic game are the anti-coordination type, i.e., its incentive to social distance increases with more of its neighbors flouting protocols and hence can be captured using the 
%linear threshold type \red{\cite{kempe2003maximizing}}
utility function in \eqref{eq_utility}. 
%\blue{CE: I am confused by this citation in the middle of the sentence. Is this a paper that has a bilinear payoff?}
The learning constants $c_i, i\in V$ represent the sensitivity of the agents to the neighbor influence.  
%\blue{CE: Avoid adverbs. They do not add meaning. Reduces clarity. } Sure!

% \red{Players want to avoid disease transmission \cite{eksin2017disease}. Each player is either healthy ($s_i = 0$) or sick ($s_i = 1$). The network $\ccalG_B$ is a contact network with each edge representing a chance of disease transmission between a healthy and a sick player. The action space captures the social distancing level of a player with action $a_i = 0$ representing self-isolation and action $a_i = 1$ representing resuming normal activity. Actions between 0 and 1 represent different levels of disease prevention measures, e.g., covering cough, or washing hands often. Resuming normal activity is the preferred action. However, if both players at the two ends of an edge take action 1, then there is a chance of disease transmission. Accordingly, the constant $c_0$ captures a healthy player's sensitivity for avoiding a risky interaction. The constant $c_1$ captures a sick player's sensitivity to avoid transmitting the disease to one of its healthy neighbors. } \blue{Shorten/Modify according to this paper}
\end{example}

\begin{example}[Political polarization] \label{example2}
The network represents the social interactions among players in opposing beliefs (agents on different sides of the bipartite network) that want to differentiate their actions from those with opposing beliefs \cite{mcconnell2018economic}. Action 1 represents a monetary choice or support for a cause that is individually desirable in the absence of partisanship. A player's tendency to take the preferred action (action 1) reduces as it has more neighbors that take action 1. That is, a player can opt-out from individual benefits or societal impact to express partisan preferences. The payoff constants $c_i$ capture the inclination of players to distinguish their actions from those in the opposite camp.

\end{example}
%Anti-coordination in social interactions naturally models a variety of scenarios other than the instances we described above 
%\green{CE: No need for perfect tense in general except in the introduction. Use present or past only in main text}. 
%\red{For a list of more examples, see \cite{bramoulle2007anti}.} \blue{CE: This citation is fit for the introduction or when you first talk about network games. Here you can cite a paper that considers network games with strategic substitutes (if such a paper exists) or one that consider scenarios relevant to the game setting or eq. (1).} 

% \begin{example}[Hawk-Dove network game] Two competing species ($\ccalS_0$ and $\ccalS_1$) face-off in an ecological environment. At each interaction players decide to be hawkish ($a_i =1$) or dovish ($a_i = 0$). A hawk move gets the highest reward if its neighboring competitors play dove. If both interacting players play dove, they miss the opportunity to overcome their competitor. If both interacting players are hawkish, they challenge each other and face costs. The constants $c_0$ and $c_1$ represent the costs species 0 and 1 incur, respectively, when they act hawkish against a hawkish competitor. 
% \end{example}

% %%%%%%%%%%%%%%%%%%%%%%%%%%%%%%%%%%%%%%%%%%%%%%%%%%%
%%%%%%%%%%%%%%%% N E W     S E C T I O N %%%%%%%%%%%%
%%%%%%%%%%%%%%%%%%%%%%%%%%%%%%%%%%%%%%%%%%%%%%%%%%%%%
\section{Local learning dynamics} \label{sec_learn}

We consider decentralized learning dynamics based on the notion of iterated elimination of dominated strategies \cite{menache2011network}. At each stage $t=1,2, \dots$, we assume agents observe the past actions of their neighbors $a_{n(i)}^{t-1}$, and determine its action $a_i^t$ according to the following rule
\begin{equation}\label{eq_updates}
\begin{aligned}
    a_i^t = 1 &\textit{, if }\qquad {1}=BR_i(\lceil a_{n(i)}^{t-1} \rceil),\\
    a_i^t = 0 &\textit{, if }\qquad {0}=BR_i(\lfloor a_{n(i)}^{t-1} \rfloor),\\
    a_i^t = \epsilon &\textrm{, otherwise } 
\end{aligned} 
\end{equation}
where $BR_i(a_{n(i)}):= \argmax_{a_i \in [0,1]} u_i(a_i,a_{n(i)})$ is the best response action profile, and $\epsilon \in (0,1)$ is an arbitrary action between 0 and 1. In \eqref{eq_updates}, agent $i$ respectively evaluates the best response function given an overestimate (ceil) and an underestimate (floor) of the sum of its neighbors' actions. The best response action for the utility function in \eqref{eq_utility} is given by
\begin{equation}\label{eq_BR}
    BR_i(a_{n(i)}) = \mathbb{1}\big(1>c_i\sum_{j\in n(i)}a_j\big) 
\end{equation}
where $\mathbb{1}(\cdot)$ is the indicator function. Accordingly, if  the overestimate of neighbors' actions is less than 1 ($c_i\sum_{j\in n(i)}\lceil a_j \rceil<1$), then agent $i$ would take the preferred action $1$ regardless of its neighbors' future actions. Similarly, if  the underestimate of neighbors' actions exceeds 1 ($c_i\sum_{j\in n(i)}\lfloor a_j \rfloor >1$), then agent $i$ would take action $0$ regardless of its neighbors future actions. In the former scenario, all actions are dominated by action 1, whereas in the latter scenario all actions are dominated by action 0. If neither of these conditions hold, i.e., when $c_i\sum_{j\in n(i)}\lfloor a_j \rfloor <1<c_i\sum_{j\in n(i)}\lceil a_j \rceil$, then agent $i$ cannot rule out any of the actions, thus it takes an arbitrary action $\epsilon\in (0,1)$. See Figure \ref{fig:bipartite1} for an illustration. 

In \cite{eksin2020control}, we provide a \textit{finite time convergence guarantee} for the dynamics discussed. Specifically, we show that the dynamics in \eqref{eq_updates} converges in at most $|V|$ iterations, eliminating all strictly dominated actions for the network game $\Gamma$ when the network $\ccalG$ is bipartite and all agents play $\epsilon$ initially, i.e. $a_i^0=\epsilon$ for $i\in V$. We note that a strictly dominated action cannot be a rational action. These updates converge to a Nash equilibrium if the game is dominance solvable, i.e., if the game is such that a single action profile is left as a result of iterated elimination of dominated strategies. For instance, the anti-coordination game $\Gamma$ with utility function in \eqref{eq_utility} is dominance solvable given constants $c_i< \frac{1}{|n(i)|}$ for all $i\in V$. Indeed, all agents take action 1 after the first update in \eqref{eq_updates} because $c_i\sum_{j\in n(i)}\lceil a_j \rceil<1$ for all $i\in V$. For general payoff constants, the game $\Gamma$ is not dominance solvable, i.e. some agents can continue to take action $\epsilon$ at the end of $|V|$ iterations. 
\section{Maximum Anti-Coordination Problem} \label{sec_MAC}

We define an edge between agents $i$ and $j$ ($(i,j)\in E$) to be inactive when at least one of the agents take action 0, i.e., when $a_i a_j=0$. Our goal is to maximize the number of inactive edges by controlling a subset of the players (with set cardinality $r\in \mathbb{Z}^+$) to play action $0$ during the learning dynamics in \eqref{eq_updates}. We state this goal to maximize anti-coordination (MAC) as follows
\begin{equation}\label{eq_MAC}
\begin{aligned} 
\max_{X\subseteq V} &\; f(X) :=\sum_{(i,j)\in E} \mathbb{1}(a_i^\infty a_j^\infty = 0) \\
\textrm{subject} \, \textrm{to} \quad &\vert X \vert = r,  \\
&a^0_j=\epsilon \quad \forall j\in V,\\
&(a^0,a^1,\ldots,a^\infty)=\Phi(a^0,X),
\end{aligned}
\end{equation}
where $\Phi(a^0,X)$ represents the sequence of actions obtained when uncontrolled agents ($V\setminus X$) follow the learning dynamics in \eqref{eq_updates}, and the actions of controlled agents are set to 0, i.e., $a_i^t=0$ for all $t\geq 0$ and $i\in X$. By removing the agents that are controlled from the game, we can guarantee that the learning process converges in finite time as per the aforementioned convergence result in \cite{eksin2020control} for bipartite networks. The control budget for the planner is restricted, i.e., the planner can only control a given $r \in \mathbb{Z}^+$ number of agents as indicated by the first constraint in \eqref{eq_MAC}.

% or in two steps for line and ring networks as per Corollary \ref{cor_convergence}. \\
In the context of disease spread in a population (Example \ref{example1}), maximizing anti-coordination in the underlying relationship network by inactivating disease transmission links is  desirable as an effective means of curbing spread of the disease between members of society. The decentralized learning dynamics do not inactivate all edges on convergence and thereby a central planner would need to control/enforce certain agents to coordinate with policy guidelines (playing action level $0$) so that the maximum number of transmission pathways shall be dismissed. 
% \blue{CE: Now motivate MAC in the context of the examples we described in the previous section.   }
% \red{\red{Agents in set $V_0$ are sick, while the others are healthy. The actions available to the agents $a_i \in [0,1]$ evolve according to the local learning dynamics established in (\ref{eq:LL}) and (\ref{eq:BR}). On convergence of the said dynamics, (we are guaranteed to converge, see Eksin \textit{et al}. \cite{eksin2018control}) we reach an action profile which need not have all disease transmission links between agents inactivated. }\blue{CE: This kind of discussion should follow after you present the problem} 

%
% REMOVE THEOREM 1
% NP-Hardness of MAC
%
\begin{theorem}\label{thm_hard}
The MAC problem in \eqref{eq_MAC} is \textrm{NP}-hard for general graphs.
\end{theorem}
\begin{proof}
We consider the NP-Complete problem of whether a vertex cover $C$ of a general graph $\ccalG$ with cardinality $|C| = k$ exists or not\begin{footnote}{A vertex cover $C$ of a graph $\ccalG=(V,E)$ is such that every edge has at least one endpoint incident (vertex) in $C$---see \cite{williamson2011design} for a formal definition.}\end{footnote}. Given a vertex cover problem $\ccalG(V,E)$ and cardinality $k$, construct a MAC for $\ccalG$ where $c_i<d_i^{-1}$ where $d_i$ is the degree of node $i\in V$, and $r=k$. Given the payoff constants $c_i$, all agents play 1 given the uncontrolled dynamics, i.e., $a_i^t = \bbone$ for all $t\geq 1$. If there exists a solution to MAC where the objective value $g(X)=|E|$, then $C=X$ is a vertex cover with cardinality $k$. If there does not exist a solution to MAC where the objective value $g(X)=|E|$, then we can conclude that there does not exist a solution to the vertex cover problem for the graph $\ccalG$ with $k$ vertices. Thus if we can solve MAC efficiently, we can solve every vertex cover problem efficiently which is a contradiction. 
% Let us consider the decision version of the Minimum Vertex-Cover problem, which is known to be NP-complete. 
% \begin{definition}(VERTEX COVER)
% Given a graph G=(V,E), does there exist a vertex cover $VC$ of $k$ vertices?
% \end{definition}
% \blue{CE: Vertex Cover is NP complete for general graphs but not for bipartite graphs. Make sure you show that your construction below is not a bipartite graph. Because the original graph is bipartite. }
% Let us say that a vertex cover for the graph $G$ of size $k$ exists. We create now a MAC instance with graph $G^{'}=G$, the learning constants chosen arbitrarily, and  cardinality of control set being fixed to $k$. The solution for this MAC instance is $\vert E \vert = m$ where $E$ is the set of edges of $G^{'}$. That is because if we control the elements in the vertex cover $VC$, all edges are inactivated, i.e. $g(VC)=m$.\\
% Conversely, let us assume that no vertex cover of size $k$ exists in the graph $G$. Then MAC cannot have a solution value $m$. Let $X$ be the control set of size $k$ for which we calculate $g(X)$. Now, since no vertex cover of size $k$ exists, we can find an edge $(u,v)$ such that no end of the edge lies in $X$. Now assign learning constants to these vertices such that $c_i < 1/d_i$ where $d_i$ is the degree of node $i$ for $i\in\{u,v\}$ which makes the nodes insensitive to peer influence and thereby they end up playing action level $1$ always. This ensures edge $(u,v)$ cannot be inactivated as a result of the controlled learning dynamics.
\end{proof}

Given that MAC is NP-Hard, we provide performance guarantees for  computationally tractable greedy approaches to solving MAC for dense bipartite networks.

\section{Average submodularity of MAC in dense bipartite networks}
In a greedy approach, we obtain a solution to a cardinality constrained maximization problem $\max_{X\subseteq V, |X|\leq r} f(X)$ by selecting one element at a time, i.e., 
\begin{align} 
    u&=\argmin_{w\in V} f(G_j \cup \{ w \})  \label{eq:Greedy}\\
    G_{j}&=G_{j-1} \cup \{u\},\; \text{for}\; 1\leq j\leq r\nonumber  
\end{align}
where $G_0 = \emptyset$. The greedy approach is a computationally tractable way to build a set of maximum cardinality $r$ for solving MAC in \eqref{eq_MAC} when, in addition, we consider the learning dynamics. If the MAC is submodular, then the greedy approach can obtain a solution comparable to the optimal set $X^*$. In the following, we provide preliminary definitions of submodularity and monotonicity, and then characterize the optimality loss when we implement a greedy selection \eqref{eq:Greedy}. 
\begin{definition}
A set function $f(\cdot)$ is submodular if for any subsets $X\subseteq Y \subseteq V$, where $V$ is the ground set, and any $u \in V\setminus Y$ we have 
\begin{align} \label{eq_sub1}
    f(X\cup \{u\})-f(X) \geq f(Y \cup \{u\})-f(Y).
\end{align}
Alternatively, for any two sets $S,T \subseteq V$, $f(\cdot)$ is submodular if
\begin{align} \label{eq_sub2}
    f(S) + f(T) \geq f(S \cup T) + f(S \cap T).
\end{align}
\end{definition}
Equations \eqref{eq_sub1} and \eqref{eq_sub2} are equivalent representations of the property, and we shall use both representations in our analyses going forward. Our results hold for a dense network -- a network where the degrees of all nodes are very large.    
\begin{definition}
For a dense network, $\frac{1}{d_i}\rightarrow 0 $ $\forall i \in V$, where $d_i:=|n(i)|$ is the degree of agent $i$.
\end{definition}
The definition of a dense network is broad and captures a variety of network configurations. A bipartite network with reasonably fewer number of nodes can be dense if the probability of edge formation between nodes is high (there is a high chance that any two agents are connected). On the other hand, a sparse bipartite graph, with a low probability of edge formation between agents can still be dense if the graph has lots of nodes.
Next we define a scenario where every agent remains undecided according to the learning dynamics \eqref{eq_updates}.

\begin{definition}
A MAC-instance defined on a bipartite graph is well-behaved if the learning constants $c_i, i\in V$ satisfy $c_i \geq \frac{1}{d_i}$ where $d_i:=|n(i)|$ is the degree of agent $i$.
\end{definition}
% \green{Changed the definition of well behaved to include equality at $c_i \geq \frac{1}{|n(i)|}$, need to change some paranthesis to brackets downstream. No major effects.}
According to the definition, there is no player that prefers to take action 1 regardless of their neighbors' actions in the well-behaved scenario. See Figure \ref{fig:bipartite1} for an example where initially without any control each player continues to take action $\epsilon$. 

\begin{lemma}
Dynamics \eqref{eq_updates} converge in one-step for a well-behaved instance when no agent is controlled.
\end{lemma}
\begin{proof}
For a well-behaved instance we converge in one step when there is no control. Agents start undecided (play $\epsilon$ initially). For all $i\in V$, $1 \leq c_i\sum_{j\in n(i)} \lceil a_j \rceil = c_i |n(i)| $ is satisfied, whereas $c_i\sum_{j\in n(i)} \lfloor a_j \rfloor=0$ as everyone plays $\epsilon$ in the zeroth step. Following the rules in \eqref{eq_updates}, agents stay undecided and the dynamics converge. 
\end{proof}

\begin{remark}
Since $\frac{1}{d_i}\rightarrow 0$, $c_i\geq \frac{1}{d_i}$ holds almost surely, i.e., MAC instances defined on dense networks are well behaved. 
\end{remark}

% Reinforcement idea. 0's push 1's which push 0's
Now consider a well-behaved instance of MAC defined on a bipartite graph $G=(V,E)$ where $V=S_0 \cup S_1, S_0 \cap S_1 = \emptyset, E \subseteq S_0 \times S_1 $. When $t=0$, all agents are undecided in the absence of control and the dynamics have converged. 
Let $X \subset S_0$ be the set of agents we control initially at time $t=1$ to play action level $0$. Initiating the dynamics with the control set allows more agents to decide to play action level $0$ or $1$. We use the notation $X_t$ to represent the set of agents choosing to play action level $0$ at time $t$. Similarly, let $X'_t$ represent the set of agents playing action level $1$ at time $t$. We thus have $X_1=X$, $X'_1=\emptyset$. 
\begin{lemma} \label{lem_reinforcement}
Undecided agents in sets $S_0$ and $S_1$ can update their strategies only in alternate time steps if controlled set is a subset of one part of the bipartite graph $X\subseteq S_0$. On convergence, a subset of agents in sets $S_0$ and $S_1$ choose to play $0$ and $1$ respectively, everyone else stays undecided.
\end{lemma}
\begin{proof}
Because of the anti-coordinating effect of the dynamics in \eqref{eq_updates}, the zero-set $X_1$ encourages agents in its neighborhood to choose action level $1$. Since $X_1=X \subset S_0$, neighbors of set $X_1$ lie exclusively in $S_1$. As a result at the next time step, we have agents in $S_1$ potentially adopting action level $1$, i.e. $S_1\supseteq X'_2 \supseteq X'_1=\emptyset$. On the other hand, $X_2=X_1=X \subset S_0$, since at $t=2$ no new agents choose to play $0$. This is because at $t=1$, all zero-playing agents are in $S_0$ (control set $X=X_0$) while all agents in $S_1$ are undecided. At $t=2$ agents in set $S_0$ update against action information of agents in set $S_1$ at $t=1$, and all agents in $S_1$ are undecided at $t=1$.

Note that because we only control agents in $S_0$, the agents choosing to play $0$ and $1$ are respectively concentrated in sets $S_0$ and $S_1$. Formally, this means that $X_t \subseteq S_0$ and $X'_t \subseteq S_1$ for all $t$. At time $t=3$, agents in $S_0$ update against the action information of agents in $S_1$ at $t=2$, and since potentially $X'_2\neq \emptyset$ (we have new agents choosing to play $1$), dynamics \eqref{eq_updates} enable $X_3 \supseteq X_2$, $X_3\subseteq S_0$. This trend will continue. At time $t=4$, $S_0\supseteq X_4=X_3$ while  $S_1 \supseteq X'_4 \supseteq X'_3$ and so on. 
Following this reasoning, we see that undecided agents in sets $S_0$ and $S_1$ can update their strategies only in \textit{alternate} time-steps. On convergence, we have $i\in X_\infty$ playing $0$,  $i\in X'_\infty$ playing $1$, rest undecided, i.e. agents in $S_0$ and $S_1$ either play $0$ or $1$ respectively, or stay undecided. Also,
\begin{align}
    X_t \subseteq X_{t+1} \subseteq \cdots \subseteq X_\infty \subseteq S_0 \nonumber \\
    X'_t \subseteq X'_{t+1} \subseteq \cdots \subseteq X'_\infty \subseteq S_1.
\end{align} \label{eq_reinforcement}
\end{proof}

This unique nature of the anti-coordination effect in well-behaved instances of MAC arises as a response to the asymmetric control that we exert, i.e. enforce $X\subset S_0$ to play $0$ at $t=1$. In the context of the epidemic game, controlling one type of agent can mean that we can only choose to isolate individuals who are sick. The agents playing zero on one side of the bipartite graph (set $S_0$) reinforce the decisions of agents playing $1$ on the other side of the bipartite network (set $S_1$), and vice versa, till there is convergence of the game. We call this a \textit{decision cascade}. Having learning constants in the range $c_i \in [\frac{1}{d_i},1)$ for all $i\in V$ allows for the removal of \textit{extremely insensitive} agents  ($c_i < \frac{1}{d_i}$) who choose to play $1$ no matter what their neighbors are doing. In the absence of these agents, the only \textit{decision cascade} is triggered by the asymmetric control we exert at $t=1$. 
% Moreover, recall from proof of Corollary \ref{cor_convergence} that there can be no agents choosing to play $0$ in the first step. 
This is the foundation of the well-behaved scenario, where agents playing $0$ and $1$ are concentrated on either side of the bipartite network.

For an illustration, see Figure \ref{fig:bipartite1} where the dynamics unfold on a well-behaved instance for MAC and the control set is confined to one side of the network. See that agents playing strategy levels $0$ and $1$ are concentrated on either side of the network. Initially at $t=1$ all agents are undecided, which is a convergent action profile by itself in the absence of control. We introduce control in the next time step which triggers the dynamics and leads to agents updating their strategies to reach a new equilibrium profile for the dynamics.

% , a property that makes anti-coordination almost akin to two parallel diffusion processes, one for agents playing $0$, one for agents playing $1$ on either side of the bipartite network, each reinforcing the other, each side making the other side more polarized and decided.  \blue{CE: split into two or more sentences. Some of the comments can be removed. You may consider having similar sentences in the introduction. Here this remark is too long and distracting. }
Next we state the main results of this section.

\begin{theorem} \label{thm_monotonic}
The MAC problem is monotone increasing for well-behaved instances when control is restricted to one-side of the bipartite graph. Given selections of control sets $X\subseteq Y \subseteq S_0$ and $f(\cdot)$ defined as in \eqref{eq_MAC}, we have 
\begin{align}
    f(X) \leq f(Y).
\end{align}
\end{theorem}
\begin{proof}
Given $X \subseteq Y$, we claim that any $X_t \subseteq Y_t$, where $X_t$ refers to the set of agents in $S_0$ that choose strategy $0$ at time $t$. Initially, the property is true ($t=0$), i.e. $X_0=X \subseteq Y=Y_0$. All agents in $S_1$ are undecided, whereas agent $i \in S_0$ plays $0$ if $i\in Y$ or $i\in X$.  Assume now that the property is true for some $t$, and let $X'_t$ and $Y'_t$ represent the corresponding sets of agents playing $1$ in $S_1$ respectively. Following the dynamics, agents in set $S_1$ update at $t+1$ according to the actions of agents in $S_0$ at time $t$. Since $X_t\subseteq Y_t$, then $X^{'}_{t+1} \subseteq Y^{'}_{t+1}$ as any agent $i\in S_1$ which plays $1$ under the influence of zero-set $X_t$ will have to choose $1$ under the influence of zero-set $Y_t$. 
This again implies that $X_{t+2} \subseteq Y_{t+2}$.
Since $X_{t+1} = X_t$ and $Y_{t+1}=Y_t$ (as the sets $S_0$ and $S_1$ update \textit{alternately}) the property holds for all $t$. Eventually, we have $X_\infty \subseteq Y_\infty$ on convergence. Since all agents playing strategy $0$ are in set $S_0$ and agents in $S_1$ are either $1$ or undecided, the only edges deactivated are the ones incident to agents in set $X_\infty$ and $Y_\infty$. Thus $f(X) \leq f(Y)$.
\end{proof}
\begin{corollary} \label{cor_monotonic}
The set of zeros $X_t$ and $Y_t$ for the processes started with zero sets $X\subseteq Y\subseteq S_0$ satisfy
\begin{align}
    X_t \subseteq Y_t \subseteq S_0.
\end{align}
\end{corollary}
\begin{theorem}\label{prop_avgsubmodular}
For a MAC instance defined on a dense bipartite network, the problem is submodular in expectation almost surely, provided control is restricted to one side of the bipartite graph. Formally, given selections of control sets $A,B \subseteq S_0$, and $f(\cdot)$ defined as in \eqref{eq_MAC} we have
\begin{align}
    \mathbb{E}[f(A)+f(B)] \geq \mathbb{E}[f(A \cup B) + f(A \cap B)] \label{subexp}
\end{align}
where the expectation is over the random learning constants $c_i$ for all $i\in V$.
\end{theorem} 
% Talk about the influence function now

The rest of the section is devoted to proving the above theorem. The proof relies on first showing that the one-step influence on an undecided agent is submodular (see Definition \ref{def_influence} and Proposition \ref{thm_submodularinfluence}) and then using a coupling argument that leverages the submodularity of the one-step influence function (see Section \ref{sec_proof_theorem}). 

We begin by defining the one-step influence function. 
% Assume we have a well-behaved instance where all agents are undecided initially. WLOG we enforce some agents in set $S_0$ to play $0$. The dynamics now unfold, and undecided agents in sets $S_0$ and $S_1$ can update in \textit{alternate} time steps as we show earlier.
\begin{definition} \label{def_influence}
Given the dynamics in \eqref{eq_updates}, let 
\begin{equation}
   f_{v}(S)= \sum_{j\in n(v)} \lfloor a_{j}(t+1)\rfloor
\end{equation}
be the one-step influence experienced by an undecided agent $v \in S_0, v\not\in S$ at time step $t+2$, where $S \subset S_0$ is the zero-set at time $t$.
\end{definition}
If the net influence $c_v f_v(S)$ on agent $v$ exceeds $1$, then agent $v$ switches to playing $0$, i.e. $a_v(t+2)=0$ following the dynamics in \eqref{eq_updates}. Next, some properties of the one-step influence function that will be useful in the proof of Theorem \ref{prop_avgsubmodular} follow. 
\begin{lemma} \label{lem_closedforminfluence}
The influence function is given by 
\begin{align}
f_v(S)=g^{v} \mathbb{1}(\bbone_n - diag({\bf c})\Mat{G}\lceil a^{S} \rceil)    
\end{align}
% \blue{CE: Very odd to use $\theta$, for the indicator function. I'd use $\mathbb{1}$ or $\bbone$. }
where $g^v$ is the row vector corresponding to the $v$'th row of the adjacency matrix $\Mat{G}$ of the underlying graph $G$,
% \red{$diag(\textbf{c})$ is a diagonal matrix with $diag(\textbf{c})_{ii}=c_i$ (the learning constant for agent $i\in V$)} 
diag(\textbf{c}) is a diagonal matrix where $\textbf{c}=[c_1,...,c_n]$ represent the learning constants for agents $i\in V$
% \blue{CE: For the diagonal matrix you can use $\diag({\bf c})$ where ${\bf c} = [c_1,\dots c_{n}]$}, \red{$1_N$} \blue{CE: Use $\bbone_N$} 
is a vector of $1$'s in $\mathbb{R}^n$, and 
% \red{$\theta(\cdot)$ is the indicator function where $\theta(x)=1$ if $x>0$, $\theta(x)=0$ otherwise}. 
$\mathbb{1}(\cdot)$ is the indicator function where $\mathbb{1}(x)=1$ if $x>0$, $\mathbb{1}(x)=0$ otherwise.
The vector $a^S \in \mathbb{R}^n$ has $a^S_i=0$ if $i \in S$, $a^S_i=\epsilon$ otherwise. 
\end{lemma}
The proof is in the Appendix.
Lemma \ref{lem_closedforminfluence} allows us to write down the closed form equation for the one-step influence exerted on an undecided agent in set $S_0$ which may (or may not) convince the agent to choose action level $0$. Writing out the influence exerted on an agent facilitates proof of Proposition \ref{thm_submodularinfluence}.
\begin{proposition} \label{thm_submodularinfluence}
The influence function $f_v(\cdot)$ is almost surely submodular in dense networks.
\end{proposition}
The proof is in the Appendix.
Submodularity of the one-step influence function eventually enables submodularity of the function $f(\cdot)$ described in \eqref{eq_MAC}. 

Next we develop the mathematical machinery needed for the coupling argument. When we initiate a \textit{decision cascade} by controlling certain agents (let's say $S\subset S_0$), the randomness of the learning constants imparts stochasticity to the set of agents that choose to play zero over time. 
% \red{Consider first the following lemma.} \blue{This lemma reads like a technical lemma. }

% \blue{I'd move lemma 4 to the appendix. }
% \red{Our objective, therefore, is straightforward. If we can construct a coupling with the properties detailed in Lemma \ref{lem_process}, we have essentially shown Theorem \ref{prop_avgsubmodular} to be true. In the rest of the paper we will endeavor to do just that. }  

The set of agents playing $0$ in $S_0$ updates itself in alternate time steps in well-behaved instances as per Lemma \ref{lem_reinforcement}. Therefore, we consider  $k=\lfloor t/2 \rfloor$ as the \textit{update time}, where $t$ represents time in \eqref{eq_updates}. Since the dynamics converge in $n$ time steps, it is enough to consider the time interval $0\leq k \leq n-1$. 
If we run the dynamics in \eqref{eq_updates} given a zero-set $S\subset S_0$, we represent the distribution of the stochastic process $\mathbb{T}=(T_k)_{k=0}^{n-1}$ where $T_k \subset S_0$ represents the zero-set at time $k$, $T_0=S$ with $\ccalQ(S)$.

% We use the notation $\mathcal{Q}(S)$ to represent the distribution of the stochastic process that determines the random set of agents playing action level $0$ over time. The process is triggered by the asymmetric control $S\subset S_0$ we apply at $t=0$.
% Given a zero-set $S\subset S_0$, if we run the dynamics in \eqref{eq_updates} then we use the notation $\mathcal{Q}(S)$ to represent the distribution of the stochastic process $\mathbb{T}=(T_k)_{k=0}^{n-1}$ where $T_k \subset S_0$ represents the zero-set at time $k$, $T_0=S$. 

In the following, we prove equivalence of the process when we consider $R^{(1)},\ldots,R^{(K)}$ a partition of $S$, i.e. $\bigcup_{i=1}^{K} R^{(i)}=S$, and sequentially apply control to the sets $R^{(j)}$, starting from $R^{(1)}$ by letting the process run with control set till convergence for each $j=2,\ldots,K$ before adding $R^{(j)}$ to the control set.
% Similarly if $R^{(1)},\ldots,R^{(K)}$ is a partition of $S$, i.e. $\bigcup_{i=1}^{K} R^{(i)}=S$, consider sequentially applying control to the sets $R^{(j)}$, starting from $R^{(1)}$ by letting the process run with control set $\bigcup_{i=1}^{j-1} R^{(i)}$ till convergence for each $j=2,\ldots,K$ before adding $R^{(j)}$ to the control set. 
Let $\mathbb{V}=(V_k)_{k=0}^{Kn-1}$ represent this process where $V_k \subset S_0$ represents the zero-set at \textit{update time} $k$. 
% More formally, we select $c_v$ for each $v\in V$ uniformly and independently in the range $(\frac{1}{d_v},1)$ (since we have a well-behaved instance). 
Then, for $1\leq j \leq K$, we set $(V_k)_{k=(j-1)n}^{jn-1} \sim \mathcal{Q}(V_{(j-1)n-1}\cup R^{(j)})$, where $V_{-1}=\emptyset$. We denote the distribution of $\mathbb{V}$ by  $\mathcal{Q}(R^{(1)},\dots,R^{(K)})$. Our first step is to show that processes $\mathbb{T}$ and $\mathbb{V}$ lead to the same end result on convergence. 
% \red{these two processes}\blue{CE: Which two processes?} lead to the same end result on convergence.
\begin{lemma} \label{lem_piecemeal}
$V_{Kn-1}$ has the same distribution as $T_{n-1}$ for any partition of the zero set $\bigcup_{i=1}^{K} R^{(i)}=S\subset S_0$.
\end{lemma}
\begin{proof}
We have $\mathbb{T}=(T_k)_{k=0}^{n-1} \sim \mathcal{Q}(S)$ and $\mathbb{V}=(V_k)_{k=0}^{Kn-1} \sim \mathcal{
Q}(R^{(1)},\cdots,R^{(K)})$. Let $\mathbb{V}'=(V'_k)_{k=0}^{Kn-1} \sim \mathcal{Q}(S,\emptyset,\cdots,\emptyset)$. Let $\mathbb{V}''=(V''_k)_{k=0}^{Kn-1} \sim \mathcal{Q}(\emptyset,\cdots,\emptyset,S)$. Since the process is monotonic (see Theorem \ref{thm_monotonic}, Corollary \ref{cor_monotonic}), by induction on the $K$ stages we have
\begin{align}
    V''_{Kn-1} \subseteq V_{Kn-1} \subseteq V'_{Kn-1}
\end{align}
However, we have $V''_{Kn-1}=V'_{Kn-1}=T_{n-1}$. Thus we have $T_{n-1}=V_{Kn-1}$.
\end{proof}
% We get the same final distribution of agents playing action level $0$ irrespective of whether we start out with the entire control set $S$ or work through its partitions one by one in stages. Consider now the following definition.
According to the above result, it does not matter whether we control all the agents in our control set right from the start, or we apply control progressively in stages by working through partitions, as eventually either process will lead to the same distribution of the agents playing action-level $0$ on convergence. Next, we consider $S\cup T\subset S_0$ as the initial zero-set. As per the above result, we can consider a partition of the set $S$, $R^{(1)},...,R^{(K)}$, and apply control progressively for all stages till all of $S$ is controlled, i.e. starting from $R^{(1)}$, let the dynamics run with control sets $\bigcup_{i=1}^{j-1} R^{(i)}$ before adding $R^{(j)}$ to the control set, and do this till we reach $S=\bigcup_{i=1}^{K} R^{(i)}$ as the control set. Next, we can add nodes from $T$ to obtain the same process as initially controlling $S\cup T$ again by the above result. In the following definition, we consider an alternative \textit{selection rule} instead of adding $T$ to the control set. 

% Now we shall reveal an even finer detail of this anti-coordination stochastic process. Since we can break the control set in partitions and apply control progressively without disturbing the distribution of the convergent profile of zero playing agents, consider the case where $S\cup T\subset S_0$ is our overall control set, and $R^{(1)},...,R^{(K)}$ is a partition of set $S$. Apply control progressively for all stages till all of $S$ is controlled, i.e. starting from $R^{(1)}$, let the dynamics run with control sets $\bigcup_{i=1}^{j-1} R^{(i)}$ before adding $R^{(j)}$ to the control set, and do this till we reach $S=\bigcup_{i=1}^{K} R^{(i)}$ as the overall control set. On convergence, instead of adding $T$ to the control set and running the dynamics for one more stage, we select agents to add to the zero-set using some carefully constructed \textit{selection rule}. We will argue next that the unique \textit{selection rule} we propose is distribution preserving, which means that following the rule we would reach the same distribution of the random set of agents playing action level $0$ on convergence as when we add $T$ to the pre-existing control set and running the dynamics. See the following Definition.
\begin{definition} \label{def_antisense}
Let $R^{(1)},\cdots,R^{(K)}$ be a partition of $S \subset S_0$ and let $T \subseteq S_0\setminus S$. Let $k=\lfloor t/2 \rfloor$ be the indicator for \textit{update time}. We define the process $\mathbb{T^{-}}=(T^{-}_k)_{k=0}^{(K+1)n-1} \sim \mathcal{Q}_{-}(R^{(1)},\cdots,R^{(K)};T)$ where, 
\begin{enumerate}
    % \item \red{for each $v\in V$ pick $c_v$ uniformly in $(\frac{1}{d_v},1)$} \blue{CE: Make this an assumption of Lemma 6}
    \item let $(T^{-}_k)_{k=0}^{Kn-1} \sim \mathcal{Q}(R^{(1)},\cdots,R^{(K)})$
    \item set $T^{-}_{Kn}$=$T^{-}_{Kn-1} \cup T$
    \item at time $(Kn+1)\leq k \leq (K+1)n-1$ initialize $T^{-}_k=T^{-}_{k-1}$ and add to $T^{-}_k$ the set of nodes in $S_0\setminus T^{-}_{k-1}$ such that
    \begin{align} \label{eq_selectionrule}
        \frac{1}{d_v} + \frac{1}{f_{v}(T^{-}_{Kn-1})}-\frac{1}{f_{v}(T^{-}_{k-1})} > c_v
    \end{align}
    until we run out of nodes to add.
\end{enumerate}
We refer to the left-hand side of Rule (3) as the selection quotient.
\end{definition}
% That is, we run the process initially with zero-sets $R^{(1)},\cdots,R^{(K)}$ sequentially. Finally, instead of adding $T$ to the resultant zero-set and running the dynamics for another stage, we select agents to add following the rule in Step $4$ in Definition \ref{def_antisense}. The following lemma establishes that even though we select agents according to the rule in Step $4$, the overall distribution of the agents playing $0$ in set $S_0$ is not affected by it.
We show next that the \textit{selection rule} we propose is distribution preserving, which means that following the rule we would reach the same distribution of the random set of agents playing action level $0$ on convergence as when we add $T$ to the pre-existing control set and run the dynamics. 

\begin{lemma} \label{lem_antisense}
The process $\mathbb{T^{-}} \sim \mathcal{Q}_{-}(R^{(1)},...,R^{(K)};T)$ as defined in Definition \ref{def_antisense} has the same distribution as the process $\mathbb{T}=(T_k)_{k=0}^{(K+1)n-1} \sim \mathcal{Q}(R^{(1)},\cdots,R^{(K)},T)$.
\end{lemma}
\begin{proof}
We are given the learning constant $c_v$ selected uniformly in $[\frac{1}{d_v},1)$ for each $v\in V$. We use the same learning constants for both processes.
Both $\mathbb{T^{-}}$ and $\mathbb{T}$ progress exactly the same way
% \blue{CE: What do you mean similarly?}
till time $k=Kn$. That is for all $0\leq k \leq Kn$ we have $T_k=T^{-}_k$. Now consider an agent $v \in S_0$ such that $v\not\in T_{k-1}$. An agent $v\notin T_{k-1}$ if $1\geq c_v f_v(T_{k-2})$ and $v\in T_k$ if $1 < c_v f_v(T_{k-1})$. Combining the two inequalities, we have 
\begin{align}
    \frac{1}{f_v(T_{k-1})}< c_v \leq \frac{1}{f_v(T_{k-2})}.
\end{align}
Conditioning on the fact that $v \not\in T_{k-1}$, we have $c_v$ uniformly distributed in $[\frac{1}{d_v},\frac{1}{f_v(T_{k-2})}]$. Consider 
\begin{align} \label{shadowLearner}
    c_{v}':=\frac{1}{d_v}+\frac{1}{f_v(T_{k-2})}-c_v
\end{align}
Then $c_{v}'$ is also distributed uniformly in $[\frac{1}{d_v},\frac{1}{f_v(T_{k-2})}]$. We call $c_{v}'$ the shadow learning constant for agent $v$. Now, if $c_{v}' > \frac{1}{f_v({T_{k-1}})}$ were also satisfied, then the distribution of zero-set $T_{k}$, conditioned on the knowledge of $T_{k-1}$ is not altered if instead of selecting nodes based on dynamics in \eqref{eq_updates} we select $v$ such that the shadow learning constant $c_{v}' > \frac{1}{f_v({T_{k-1}})}$. That is, instead of checking whether $c_v f_v(T_{k-1})>1$, we check whether $c_{v}'f_v(T_{k-1})>1$. This means that if we pick agent $v$ using the shadow learning constant $c_v'$ rather than its original learning constant $c_v$, we preserve the distribution of $T_k$.  Plugging in $c_{v}'$ from \eqref{shadowLearner}, we get
\begin{align} \label{eq_selrule0}
    \frac{1}{d_v}+\frac{1}{f_v(T_{k-2})}-c_v > \frac{1}{f_v(T_{k-1})}.
\end{align}
Since both processes $T$ and $T^{-}$ run the same way for the first $K$ stages, we get $T_{Kn-1}=T^{-}_{Kn-1}$. Moreover, in either case we initialize $T_{Kn}=T^{-}_{Kn}=T_{Kn-1}\cup T$. Now we follow a very similar reasoning as earlier. For an agent $v \not \in T_{Kn}$ but $v\in T_k$, for $Kn+1 \leq k \leq (K+1)n-1$, we need 
\begin{align}
    \frac{1}{f_v(T_{k-1})}< c_v \leq \frac{1}{f_v(T_{Kn-1})}
\end{align}
Conditioning on $v\not \in T_{Kn}$ we have $c_v$ uniformly distributed in $[\frac{1}{d_v},\frac{1}{f_v(T_{Kn-1})}]$.
Define the shadow learning constant
\begin{align}
    c_{v}'= \frac{1}{d_v}+ \frac{1}{f_v(T_{Kn-1})} - c_v.
\end{align} 
Then, $c_{v}'$ is also distributed uniformly in $[\frac{1}{d_v},\frac{1}{f_v(T_{Kn-1})}]$. If $c_{v}'>\frac{1}{f_v(T_{k-1})}$, then selecting agent $v \not \in T_{k-1}$ to be added to zero-set $T_k$ keeps the distribution of $T_k$ the same as if we had followed the dynamics in \eqref{eq_updates}. Plugging in the value of $c_{v}'$ we get the equivalent of \eqref{eq_selectionrule} for process $\mathbb{T}$, which is guaranteed to not alter the distribution of $T_k$ if we would just follow the dynamics \eqref{eq_updates}.
% Since the game is progressive, we use the knowledge of the set $T^{-}_{k-1}$ to additionally enforce that $v\in S_0\setminus T^{-}_{k-1}$ to recover \eqref{eq_selrule0} from \eqref{eq_selrule}.
% This does not alter the distribution of set $T^{-}_k$ from that of $T_k$ obtained from process $\mathbb{T}$.
\end{proof}

\begin{remark}

Since the zero sets $T^{-}_{Kn-1} \subseteq T^{-}_{k-2}$, from monotonicity of $f_v(\cdot)$ we get $\frac{1}{f_v(T^{-}_{Kn-1})}\geq \frac{1}{f_v(T^{-}_{k-2})}$. Thereby the selection rule \eqref{eq_selectionrule} is an overestimate of what we need to exactly determine whether $v\in T^{-}_k$. The rule additionally prescribes such $v\not\in T^{-}_{Kn}$ which may have already been added to the zero-set before time $k$ and hence we additionally prescribe $v\not \in T^{-}$ as the dynamics are progressive.
\end{remark}

Next we provide two technical results related to monotone submodularity and coupling which we use in the proof of Theorem \ref{prop_avgsubmodular}.
\begin{lemma} \label{lem_consequence}
Let $h:2^{V}\rightarrow \mathbb{R}^{+}$ be monotone and submodular. Then if $I\subseteq I'\subseteq V$ and $J\subseteq J' \subseteq V$ are given we have 
\begin{align} \label{eq_consequence}
    h(I\cup J') -h(I) \geq h(I'\cup J)-h(I').
\end{align}
\end{lemma}
% \begin{proof}
%  See \cite{mossel2010submodularity} for a proof.
% \end{proof}
See \cite{mossel2010submodularity} for a proof.
\begin{lemma}\label{lem_process}
Consider $A,B \subseteq S_0$. Let $Z=A\cap B$ and $D=A\cup B$. Let $A_{t} \subseteq S_0$ represent the zero-set at time $t$ when we initiate the dynamics in \eqref{eq_updates} with the control set $A$. We let the dynamics in \eqref{eq_updates} unfold on a MAC instance defined on a dense network with these four control sets, and then couple the following processes
\begin{align}
    \mathbb{A}= (A_t)_{t=0}^{\infty} \sim \mathcal{Q}(A), \nonumber \\
    \mathbb{B}=(B_t)_{t=0}^{\infty} \sim \mathcal{Q}(B), \nonumber \\
    \mathbb{Z}=(Z_t)_{t=0}^{\infty} \sim \mathcal{Q}(Z), \nonumber \\
    \mathbb{D}=(D_t)_{t=0}^{\infty} \sim \mathcal{Q}(D) \nonumber 
\end{align}
% \blue{CE: What does $\mathcal{Q}(\cdot)$ function represent?}
in such a way that $Z_{\infty} \subseteq A_{\infty} \cap B_{\infty}$ and $D_{\infty} \subseteq A_{\infty} \cup B_{\infty}$. If such a coupling exists for any selections of $A,B$, then we obtain  
submodularity in expectation \eqref{subexp} given the random learning constants $\textbf{c}$.
% \red{we have proven Theorem \ref{prop_avgsubmodular}.} \blue{CE: State this formally...E.g., then we have submodularity in expectation \eqref{subexp} given the random learning constants ...}
\end{lemma}
The proof is in the Appendix.
We now have all the elements in place for us to prove Theorem \ref{prop_avgsubmodular}.
% We now have all the tools in place to attack the proof of Proposition \ref{prop_avgsubmodular} head-on.
\subsection{Proof of Theorem \ref{prop_avgsubmodular}} \label{sec_proof_theorem}
Consider $A,B \subseteq S_0$. Let $Z=A\cap B$ and $D=A\cup B$. Let $A_{t} \subseteq S_0$ represent the zero-set at time $t$ when we initiate the dynamics in \eqref{eq_updates} with the control set $A$. We let the dynamics in \eqref{eq_updates} unfold on a MAC instance defined on a dense network with these four control sets, i.e. $A,B,Z$ and $D$ respectively, which lead to the following four stochastic processes recording the set of agents playing action level $0$ over time
\begin{align} \label{eq_4processes}
    \mathbb{A}= (A_t)_{t=0}^{\infty} \sim \mathcal{Q}(A), \nonumber \\
    \mathbb{B}=(B_t)_{t=0}^{\infty} \sim \mathcal{Q}(B), \nonumber \\
    \mathbb{Z}=(Z_t)_{t=0}^{\infty} \sim \mathcal{Q}(Z), \nonumber \\
    \mathbb{D}=(D_t)_{t=0}^{\infty} \sim \mathcal{Q}(D). \nonumber 
\end{align}
The four processes above can be equivalently represented as 
\begin{equation} \label{eq_representation}
\begin{aligned}
    \mathbb{A}=(A_k)_{k=0}^{3n-1} &\sim \mathcal{Q}(A\cap B,A\setminus B,\emptyset) \\ 
    \mathbb{B}=(B_k)_{k=0}^{3n-1} &\sim \mathcal{Q}_{-}(A\cap B, \emptyset; B\setminus A) \\ 
    \mathbb{Z}=(Z_k)_{k=0}^{3n-1} &\sim \mathcal{Q}(A\cap B, \emptyset,\emptyset) \\
    \mathbb{D}=(D_k)_{k=0}^{3n-1} &\sim \mathcal{Q}_{-}(A\cap B,A\setminus B; B\setminus A),
    \end{aligned}
\end{equation}
using Lemmas \ref{lem_piecemeal} and \ref{lem_antisense} and the finite time convergence guarantee for the dynamics. 
For well behaved instances, it suffices to only keep record of alternate time steps for the processes given by $k=\lfloor t/2 \rfloor$, also referred to as \textit{update time steps}. 
% \red{Our question now is whether there exists a coupling of the processes such that $Z_k \subseteq A_k\cap B_k$ and $D_k \subseteq A_k \cup B_k$ for all $0\leq k\leq 3n-1$.} \blue{CE: It seems you proposed the coupling in the above equation and with the next sentence. The goal is to show that the coupling is actually legitimate.}

We will now show that our proposed coupling in \eqref{eq_representation} will have $Z_k \subseteq A_k \cap B_k$ and $D_k \subseteq A_k \cup B_k$ for all $0\leq k \leq 3n-1$. 
We use the same learning constants for all the four processes. By construction, for all $0\leq k\leq 2n-1$,
$B_k=Z_k\subseteq A_k$ and $Z_k=A_k \cap B_k$. Similarly, for all $0\leq k \leq 2n-1$, we have $D_k=A_k$ which implies $D_k \subseteq A_k \cup B_k$. Thus for all \textit{update times} $0\leq k\leq 2n-1$ our result holds.
% \red{Well, what can we say for $2n\leq k \leq 3n-1$?} \blue{CE: State what you want to say and then we will show this by proving the following two statements.}
We need to additionally show that the inequalities hold for all times $2n\leq k \leq 3n-1$.
Consider the following statements.\\
\textit{Statement 1}: 
\begin{align}
    D_k\setminus D_{2n-1} \subseteq B_k \setminus B_{2n-1}
\end{align}
\textit{Statement 2}: 
\begin{align}
    \frac{1}{f_v(B_{2n-1})}-\frac{1}{f_v(B_k)} \geq \frac{1}{f_v(D_{2n-1})}-\frac{1}{f_v(D_k)}
\end{align}
First we shall the use the principle of mathematical induction to show that these two statements are valid for all time $2n \leq k \leq 3n-1$. 

Consider $k=2n$. For the first statement, we have $D_{2n}=D_{2n-1}\cup (B\setminus A)$, $B_{2n}=B_{2n-1}\cup (B\setminus A)$. That means we have 
\begin{align}
    D_{2n}\setminus D_{2n-1} 
    &=(D_{2n-1}\cup (B\setminus A))\setminus D_{2n-1} \nonumber\\ 
    &=\{(B\setminus A)\setminus D_{2n-1} \} \cup \{(D_{2n-1})\setminus D_{2n-1}\} \nonumber \\
    &=(B\setminus A)\setminus D_{2n-1}. \nonumber 
\end{align}
Now using $B_{2n-1}\subseteq D_{2n-1}$, we have 
\begin{align}
 D_{2n}\setminus D_{2n-1}    
    &\subseteq (B\setminus A)\setminus B_{2n-1} \nonumber \\  
    % (\because B_{2n-1}\subseteq D_{2n-1} )  \\ \nonumber
    &=\{(B\setminus A)\setminus B_{2n-1}\} \cup \{ B_{2n-1} \setminus B_{2n-1}\} \nonumber \\ 
    &= (B_{2n-1}\cup (B\setminus A))\setminus B_{2n-1} \nonumber \\
    &= B_{2n}\setminus B_{2n-1}.
\end{align}
% \blue{CE: explain each equality.}
Next we prove Statement 2 for $k=2n$. Now in \eqref{eq_consequence}, plug in $I=B_{2n-1}$, $I'=D_{2n-1}$, $J=J'=B\setminus A$. Noting that $f_v(\cdot)$ is monotone increasing (Corollary \ref{cor_monotontinfluence} in the Appendix) and submodular almost surely for dense networks (Proposition \ref{thm_submodularinfluence}), we have
\begin{align} \label{eq_iteration1}
    f_v(B_{2n})-f_v(B_{2n-1}) \geq f_v (D_{2n}) -f_v (D_{2n-1}) 
\end{align}
Equation \eqref{eq_iteration1} holds almost surely for dense networks. Also, as a consequence of the monotonicity of $f_v(\cdot)$ we have 
\begin{align} \label{eq_iteration2}
    f_v(B_{2n})f_v(B_{2n-1}) \leq f_v(D_{2n})f_v(D_{2n-1})
\end{align}
% \red{Dividing Equation \eqref{eq_iteration1} with Equation \eqref{eq_iteration2}, we get} 
Dividing the left (right) hand side of \eqref{eq_iteration1} with the left (right) hand side of \eqref{eq_iteration2}, we get
\begin{align}
    \frac{1}{f_v(B_{2n-1})} -\frac{1}{f_v(B_{2n})} \geq \frac{1}{f_v(D_{2n-1})} -\frac{1}{f_v(D_{2n})}
\end{align}
Thus statement 2 holds at $k=2n$.
Assume now that the statements hold at time $k$. Statement 2 at time $k$ is given by  
\begin{align} \label{selection_quotient}
    \frac{1}{d_v}+\frac{1}{f_v(B_{2n-1})} -\frac{1}{f_v(B_k)} \geq \frac{1}{d_v}+\frac{1}{f_v(D_{2n-1})} -\frac{1}{f_v(D_k)} 
\end{align}
Here we have added $\frac{1}{d_v}$ on either side. This inequality says that the selection quotient for the process $\mathbb{B}$ is larger than that for process $\mathbb{D}$. See that for both processes $\mathbb{B}$ and $\mathbb{D}$ we select elements using selection rule (3) in Definition \eqref{def_antisense} in the last stage (for $2n \leq k \leq 3n-1$). Since process $\mathbb{B}$ features a higher selection quotient,  
% \red{Match this with the selection rule (4) in Definition \eqref{def_antisense}.} \blue{Not sure what you mean by ``matching this''...}This implies that
any $v\in S_0\setminus D_k$ that is prescribed to be added to set $D_{k+1}$ is also prescribed to be added to set $B_{k+1}$, unless it is already in set $B_k$. Since $D_k\setminus D_{2n-1} \subseteq B_k \setminus B_{2n-1}$ by assumption, this further implies that 
\begin{align} \label{eq_stmnt1}
    D_{k+1}\setminus D_{2n-1} \subseteq B_{k+1}\setminus B_{2n-1}. 
\end{align}
% holds, since by assumption, $D_{k}\setminus D_{2n-1} \subseteq B_k\setminus B_{2n-1}$. 
This is nothing but Statement 1 holding at \textit{update time} $k+1$. We plug in $B_{2n-1} \subseteq D_{2n-1} $
 and \eqref{eq_stmnt1} in 
% \red{Lemma \eqref{lem_consequence}}
\eqref{eq_consequence} to get that 
\begin{align}
    f_v(B_{k+1}) -f_v(B_{2n-1}) \geq f_v(D_{k+1}) - f_v(D_{2n-1})
\end{align}
Also, since $f_v(B_{k+1}) f_v(B_{2n-1}) \leq f_v(D_{k+1}) f_v(D_{2n-1})$ by monotonicity of $f_v(\cdot)$, we show that Statement 2 is true at \textit{update time} $k+1$. Thus statements 1 and 2 are true for all update times, and all intermediate time steps between two update time steps. 

Now, for $2n \leq k \leq 3n-1$, we have $A_k = D_{2n-1}$. Also $D_k\setminus D_{2n-1} \subseteq B_k \setminus B_{2n-1} \subseteq B_k$ which implies that 
\begin{align}
    (D_k \setminus D_{2n-1}) \cup D_{2n-1} &\subseteq A_k \cup B_k \nonumber \\  
    D_k &\subseteq A_k \cup B_k
\end{align}
Also, $Z_k \subseteq A_k \cap B_k$ holds from the construction for this time range. As the dynamics definitely converge by $k=3n-1$ by the finite time convergence guarantee established in \cite{eksin2020control}, we have $Z_\infty \subseteq A_\infty \cap B_\infty$ and $D_\infty \subseteq A_\infty \cup B_\infty$. We have successfully constructed the coupling we needed. The  the desired relation in \eqref{subexp} follows by Lemma \ref{lem_process}.
\hfill \QED

\subsection{Performance guarantee for greedy selection}
Combining the results of Theorems \ref{thm_monotonic} and  \ref{prop_avgsubmodular}, the greedy algorithm to solve MAC will achieve the standard $1-1/e$ performance guarantee for submodular function maximization on average. %By average, here we refer to $\mathbb{E}(f(\cdot))$, where the expectation is over the random learning constants $\textbf{c}$.
%  \begin{definition}
%  By average, we measure anti-coordination as $h(X)=\mathbb{E}(f(X))$ where $f(\cdot)$ is defined as in \eqref{eq_MAC} measures the number of edges deactivated by the learning dynamics in \eqref{eq_updates} under the influence of control set $X \subset S_0$ and the expectation is over the random learning constants $\textbf{c}$. 
%  \end{definition}

 \begin{theorem} \label{thm_greedybound}
 If the control budget is $r$, then the solution obtained by the greedy algorithm after $r$ steps satisfies
 \begin{align}
    \mathbb{E}f(G_r)\geq (1-\frac{1}{e})\mathbb{E}f(X^{*})+ \frac{1}{e}(1-\frac{1}{r})\mathbb{E}f(\emptyset) 
 \end{align}
 \end{theorem}
 The proof is in the Appendix.
%  More precisely, if the control budget is $r$, then the solution obtained by the greedy algorithm after $r$ steps satisfies
%  \begin{align}\label{greedybound}
%      h(G_r)\geq (1-\frac{1}{e})h(X^{*})+ \frac{1}{e}(1-\frac{1}{r})h(\emptyset)
%  \end{align}
 The greedy bound benefits from the fact that $f(\emptyset)\not=0$, as even without any control, anti-coordinating behavior among agents following the decentralized dynamics leads to deactivation of some edges on convergence. 

% \red{By average, we mean we measure anti-coordination as the expected number of edges deactivated on convergence of the controlled learning dynamics in \eqref{eq_updates}. }\blue{Rewrite. This is a little too informal.}

% %%%%%%%%%%%%%%%%%%%%%%%%%%%%%%%%%%%%%%%%%%%%%%%%%%%
%%%%%%%%%%%%%%%% N E W     S E C T I O N %%%%%%%%%%%%
%%%%%%%%%%%%%%%%%%%%%%%%%%%%%%%%%%%%%%%%%%%%%%%%%%%%%
\begin{figure}
\centering
\includegraphics[scale=0.65]{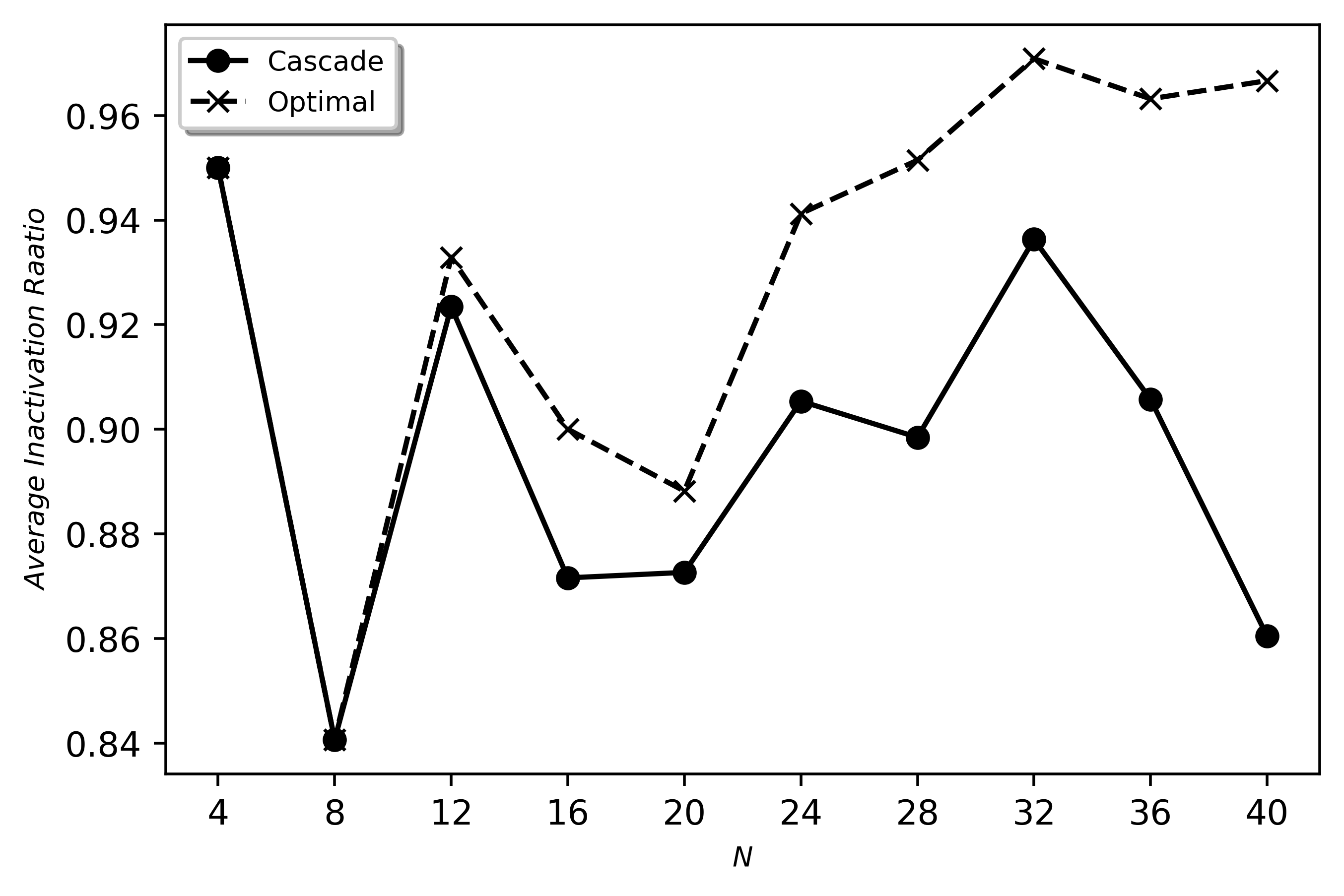}
\includegraphics[scale=0.65]{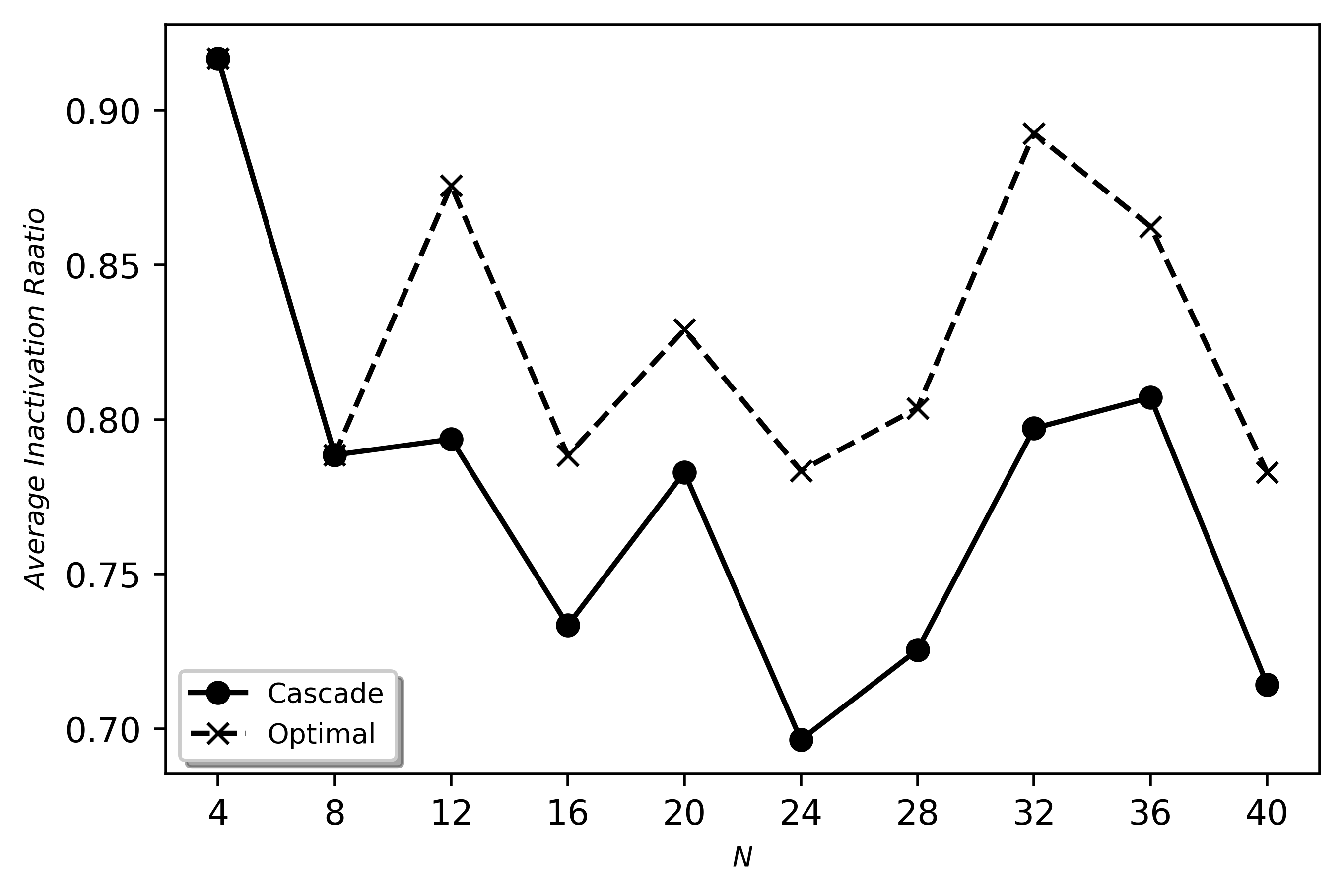}
\caption{Average Inactivation Ratio vs Cardinality ($n$) of the Bipartite Network for connection probability equal to $0.3$ and $0.8$ respectively. Control budget set at $\lceil \frac{n}{10} \rceil$ for every network realization.}
\label{fig:InactivationRatio}
\vspace{-12pt}
\end{figure}

\section{Simulation}
MAC is concerned with the deactivation of as many edges as possible on convergence of learning dynamics, perturbed by controlling a select few agents. We define the \textit{Inactivation Ratio} as the ratio of the number of edges inactivated on convergence to the number of active edges in the network before the dynamics progress. \textit{Inactivation Ratio}, therefore, is a measure of how successful MAC is on the particular graph instance, given the control.

For our simulations we consider random bipartite graph instances with edge formation probability equal to $0.3$ and $0.8$. Every realization of a network for given network sizes (\{$4,8,12,..,40$\}) has a random topology with random learning constants for the agents. We sample the learning constants $c$ for every agent from a uniform distribution between the limits $[0,1]$. The control budget is fixed at $\lceil\frac{n}{10} \rceil$, where $n=|V|$ is the number of nodes in the graph. Given the budget, we select the control profile using a greedy cascade based algorithm \eqref{eq:Greedy}. We compare its anti-coordination performance with a control set generated using brute force search. In the brute force approach, we go over all the possible control sets for the budget specified and find the one that maximizes the number of edges deactivated. For every network instance, we calculate the \textit{Inactivation Ratio} for both the control sets found using the greedy algorithm \eqref{eq:Greedy} and brute force search. For a given network size, we sample 40 instances of random bipartite graphs and evaluate the performance of the greedy algorithm. 

We plot the average \textit{Inactivation Ratio} against the size of network in Figure \ref{fig:InactivationRatio}. %We control the sparsity of the networks by varying the edge connection probability between the two plots. 
The \textit{Inactivation Ratio}, on average, for the greedy algorithm is close to the optimal inactivation at the current control budget for every network size. The maximum inactivation ratio gap for our simulations stands at $0.106$ for networks with low edge formation probability and $0.095$ for the networks with high edge formation probability, further highlighting the good performance of \eqref{eq:Greedy} in selecting control agents to induce anti-coordination. Note that in Figure \ref{fig:InactivationRatio} (Top), when the network size increases, since the connection probability is only $0.3$, we move far away from our dense network regime. As a result the gap widens between the solution we achieve using the greedy algorithm and the optimal solution. However, a larger network size with a high connection probability of $0.8$ takes us closer to the dense network regime in Figure \ref{fig:InactivationRatio} (Bottom), and the gap reduces as our performance guarantees hold.

\begin{figure}
\centering
\includegraphics[scale=0.65]{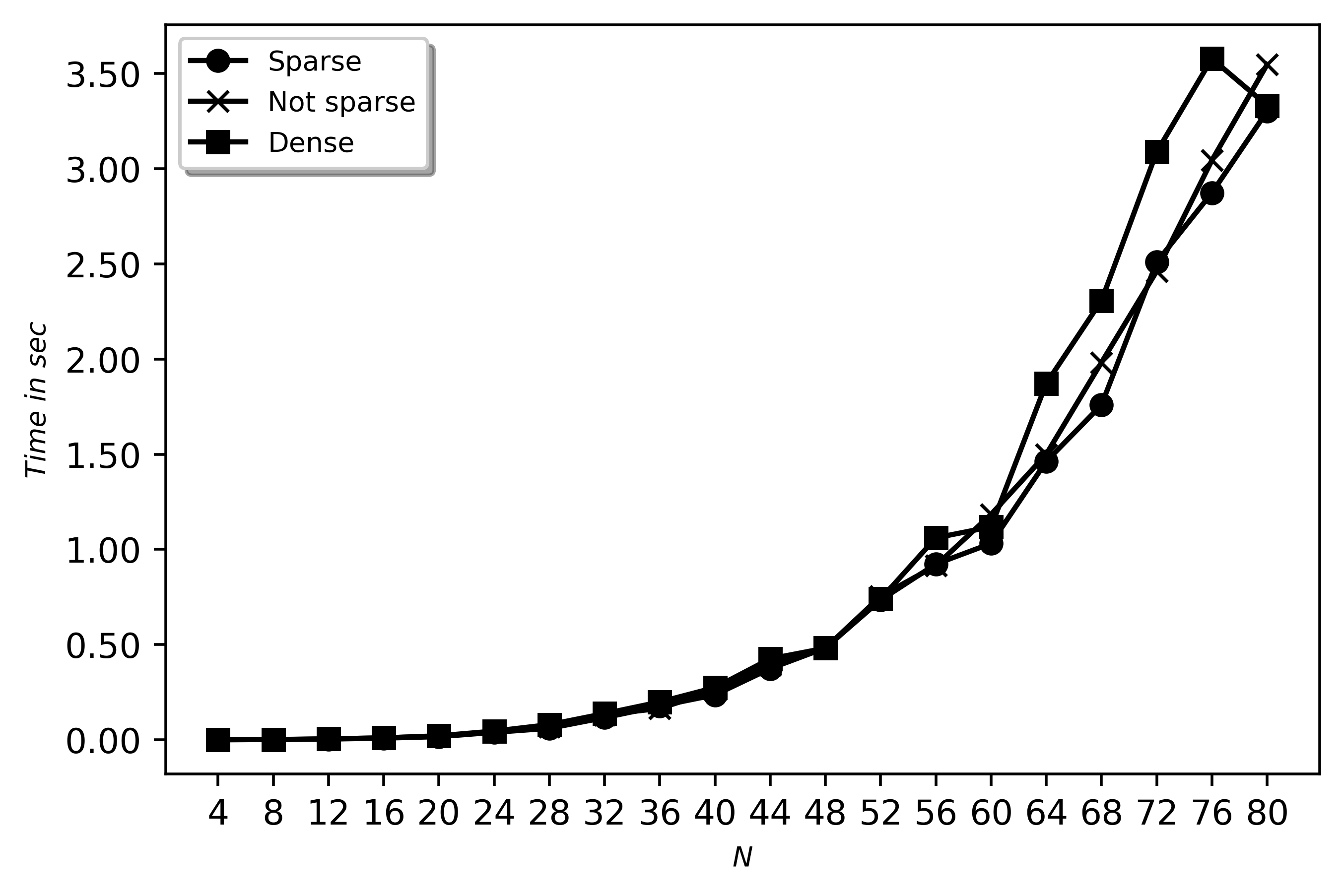}
\caption{Average time needed to construct greedy set for different network sizes and edge formation probabilities (0.3, 0.5, 0.8 for Sparse, Not Sparse and Dense Bipartite Networks respectively). Control budget fixed at $\lceil \frac{n}{10} \rceil$ where $n$ is the network size.} \vspace{-12pt}
\label{fig:TimeComplexity}
\end{figure}
In Figure \ref{fig:TimeComplexity}, we show how long it takes to build the greedy control set using \eqref{eq:Greedy} for different control budgets and network sizes. All our simulations have been performed on Apple(R) M1 CPU (Arm(R) based, 8-core) with 16GB of RAM.
% Remove
% Remove section on Computational complexity
% Remove
\subsection{Computational complexity}
 Note the greedy algorithm has a time complexity of $O(n^3)$. This is because at each step we potentially scan $O(n)$ agents for finding the one which causes the highest cascade. For each such scan, we run learning dynamics which takes at most $O(n)$ time steps to converge \cite{eksin2020control} and each step of the learning dynamics involves $O(n)$ updates. 
%  For ring and line networks, using Corollary \ref{cor_convergence}, we have a time complexity of $O(N^2)$ for the greedy algorithm.

\section{Conclusion}
%\blue{Start by defining the problem. Use past tense, e.g., We considered anti-coordination network games } 
We defined the combinatorial problem of selecting agents to control to maximize anti-coordination among rational agents in a network game. Anti-coordination is measured as the number of edges deactivated from the network on convergence of decentralized learning dynamics. Our paper is centered around establishing submodularity of 
MAC for dense bipartite networks. For that, we first establish submodularity and monotonicity of a one-step influence function. We follow that up with discussions around the stochastic nature of the dynamics, given random learning constants. We establish that the process has a unique property that instead of applying control in a consolidated fashion, we can apply it in a distributed manner which opens up possibilities for alternate equivalent descriptions of the game (Lemmas \ref{lem_piecemeal} and \ref{lem_antisense}). We then derive a selection rule which when substituted for the dynamics in the ultimate stage of control does not alter the distribution of agents playing strategy level $0$ on convergence. Using all these results in conjunction allows to derive the coupling argument which guarantees that the four processes we couple have the properties we desire in Lemma \ref{lem_process}. This leads directly to the proof of submodularity (Theorem \ref{prop_avgsubmodular}).
We also show that MAC is monotone increasing in the set of control agents.
Using these results together, we derived the approximation guarantee for greedy node selection for MAC (Theorem \ref{thm_greedybound}). Our computational results indicate that greedy selection strategies can be effective in producing near-optimal control sets for MAC on general bipartite network scenarios. %An immediate extension of our work would be to provide optimality certificates for greedy node selection in such general situations and/or quantify the worst performance we can expect to obtain.\\

% Anti-coordination appropriately models a lot of situations where decision making entities compete against each other in complex network structures. One assumption of our work is that the underlying graph is a constant. We can explore game convergence and performance of greedy node selection in situations where the interaction graph itself is time-varying. Moreover, we assume that the central planner has access to complete information of the game, the network and the utilities of the agents. It will be interesting to analyse situations where the planner has access only to partial information and/or the network structure in revealed over time such as in an online setting.

%\red{To be written}

% % % %%%%%%%%%%%%%%%%%%%%%%%%%%%%%%%%%%%%%%%%%%%%%%%%%%%
% % %%%%%%%%%%%%%%%% N E W     S E C T I O N %%%%%%%%%%%%
% % %%%%%%%%%%%%%%%%%%%%%%%%%%%%%%%%%%%%%%%%%%%%%%%%%%%%%
\section{Appendix}
\subsection{Proof of Theorem \ref{thm_greedybound}}
\begin{proof}
Let $G_j$ be the control set obtained after $j$ steps of the greedy algorithm. Define $h(\cdot):= \mathbb{E}f(\cdot)$.
Since the function $f(\cdot)$ is monotonic (Theorem \ref{thm_monotonic}) and monotonicity is preserved under expectation, we have
\begin{align}
    h(X^*) &\leq h(X^* \cup G_j)  \nonumber \\
    &= h(G_j)+ \sum_{i=1}^s [h(\tau_{i-1} \cup \{ e_i^{*} \}) - h(\tau_{i-1})]  \nonumber
\end{align}
where $\tau_0=G_j$, $\tau_i=G_j \cup \{e_1^{*},\cdots,e_i^{*}\}$ for $i=1,\ldots,s$ , $e_i^*$ is the $i^{th}$ element of $X^*$.

% and $k_j$ is the monotonicity violation at step $j$ of the greedy selection procedure. 

The second line above follows by splitting the right hand side into a telescopic sum. Now, since $f(\cdot)$ is submodular in expectation (Theorem \ref{prop_avgsubmodular}), and $G_j \subseteq \tau_i $ for all $i=1,..,s$, %we get 
\begin{equation}
    h(\tau_{i-1} \cup \{ e_i^*\}) - h(\tau_{i-1}) \leq h(G_j \cup \{ e_i^*\}) - h(G_j) 
\end{equation}
Hence,
\begin{align}
    h(X^*) &\leq h(G_j) + \sum_{i=1}^s [f(G_j \cup \{ e_i^*\})-f(G_j)]   \nonumber \\
    &\leq h(G_j) + s[h(G_{j+1})-h(G_j)]
\end{align}
where the second inequality follows by the fact that we use a greedy approach to construct $G_{j+1}$ using $G_j$ as in \eqref{eq:Greedy}. 

Now define $\delta_j \equiv h(X^*)-h(G_j)\geq 0$ as the optimality gap obtained after $j$ steps of the greedy algorithm. Then, 
\begin{align}
    \delta_j \leq s[-\delta_{j+1} + \delta_j] 
\end{align}
Moving $\delta_j$ from left to right, $\delta_{j+1}$ from right to left, and dividing both sides by $s>0$, we get the following recursion
\begin{equation} \label{eq_one_step_bound}
    \delta_{j+1} \leq (1-1/s)\delta_j 
\end{equation}
% Notice that all $\delta_j$ are nonnegative by definition. Now, writing out the terms of the recursion starting from $\delta_1$, we get 
% \begin{align}
%     \delta_1 \leq (1-1/s)\delta_0 + (\epsilon + k/s) \nonumber \\
%     \delta_2 \leq (1-1/s)\delta_1 + (\epsilon + k/s) \nonumber \\
%     \leq (1-1/s)[(1-1/s)\delta_0 + (\epsilon + k/s)] + (\epsilon + k/s) \nonumber \\
%     = (1-1/s)^2 \delta_0 + (1+(1-1/s))(\epsilon + k/s) \nonumber
% \end{align} and so on. 
Implementing the bound in \eqref{eq_one_step_bound} for the right hand side until we get to $j=0$, we get the following optimality gap bound 
\begin{align}
    \delta_r \leq (1-1/s)^r \delta_0  \label{eq_deltar}
\end{align}
Let $q=(1-1/s)$.

%and $S=\sum_{j=0}^{r-1} \frac{k_j}{s}(q)^{r-1-j}$. 
%From Proposition \ref{thm_almost_mon}, we see that 
%\begin{align}
%f(X^{*})\leq f(X^{*} \cup G_j) + 2j    
%\end{align}
%holds, assuming in the worst case that all elements in $G_j$ are not in $X^{*}$. Then, $k_j=2j$. Plugging $k_j$ in $S$, we get an Arithmetic-Geometric Progression. Using standard techniques, we yield
%\begin{align}
  %  S= 2(r-1) - 2(s-1)(1-q^{r-1}) \label{eq_S}
%\end{align}
%Again, 
%\begin{align}
   % \sum_{i=0}^{r-1}(1-1/s)^i=\ep%silon s (1-q^{r}) \label{eq_P}
%\end{align}

Substituting the values for $\delta_r$ and $\delta_0$, we get after some rearranging,
\begin{align} \label{eq_Gr_1}
    h(G_r) \geq (1-q^r)h(X^{*}) + q^{r}h(\phi) 
\end{align}     
    %-\epsilon s (1-q^{r}) \nonumber \\
    %-2(r-1) + 2sq(1-q^{r-1})

Now, since MAC is an equality constrained Optimization problem, and we run the greedy algorithm for as many steps as is the control budget $r$, we have $s=r$. We obtain, from \eqref{eq_Gr_1}
\begin{align}
    h(G_r) \geq (1-q^{r})h(X^{*}) + q^{r}h(\phi)  
\end{align}
Now, $h(X^{*}) \geq 0$. Using the inequality $q^{r}=(1-1/s)^{r} \leq e^{-r/s}$ $\forall s\geq 1$, we get
\begin{align}
(1-q^{r})h(X^{*}) \geq (1-e^{-1})h(X^{*})  
\end{align}
 Further, we can use the inequality 
$q^{r}=(1-1/r)^{r}\geq e^{-1}(1-1/r)$ $\forall r\geq 1$ to obtain  our desired result.
% (We use $(1+\frac{x}{n})^n \geq e^x(1-\frac{x^{2}}{n})$ $\forall n\geq 1, |x|\leq n$). 
% Finally, we get,
% \begin{align}
%     f(G_r) \geq (1-e^{-1})f(X^{*}) + e^{-1}(1-\frac{1}{r})(f(\phi)+2r)-4r
% \end{align}
% which is our desired result in \eqref{eq_greedy_bound}.
\end{proof}

\subsection{Proof of Lemma \ref{lem_closedforminfluence}}

% \begin{lemma}
% The influence function is given by 
% \begin{align}
% f_v(S)=g^{v}\theta(1_N - diag(\textbf{c})\Mat{G}\lceil a^{S} \rceil)    
% \end{align}
% where $g^v$ is the row vector corresponding to the $v$'th row of the adjacency matrix $\Mat{G}$ of the underlying graph $G$, $diag(\textbf{c})$ is a diagonal matrix with $diag(\textbf{c})_{ii}=c_i$ (the learning constant for agent $i\in V$), $1_N$ is a vector of $1$'s in $\mathbb{R}^N$, and $\theta(\cdot)$ is the indicator function where $\theta(x)=1$ if $x>0$, $\theta(x)=0$ otherwise. The vector $a^S \in \mathbb{R}^N$ has $a^S_i=0$ if $i \in S$, $a^S_i=\epsilon$ otherwise. 
% \end{lemma}

Given $S\subset S_0$, the zero-set at some time $t$, $f_v(S)$ measures the number of agents in $n(v)$ who are playing $1$ in time step $t+1$.  The indicator vector $\mathbb{1}(1_n - diag(\textbf{c})\Mat{G}\lceil a^{S} \rceil)_j=1$ for all neighbors $j\in n(v)$ of $v$ which play $1$ at time-step $t+1$. Also $g^v=\Mat{G}_v$, $g^v_i=1$ if $i\in n(v)$, $g^v_i=0$ otherwise. Left multiplying the indicator function $\mathbb{1}(\cdot)$  with row-vector $g^v$ ensures we count exclusively those components of $\mathbb{1}(1_n - diag(\textbf{c})\Mat{G}\lceil a^{S} \rceil)$ that correspond to the neighbors of $v$. We construct vector $a^S$ such that for all $i\in S_1$, $a^S_i=\epsilon$ and for $i\in S_0$, $a^S_i=a_i(t)$. Let $j\in n(v)$. Then $j\in S_1$, which implies that $n(j)\subseteq S_0$. Therefore, to check whether $j$ chooses strategy level $1$ at time step $t+1$, we evaluate $\mathbb{1}(1-c_j\sum_{i\in n(j)}\lceil a_i(t)\rceil)$, which corresponds to the $j$'th component of the indicator vector. Since $i\in S_0$, we do not care about the components of $a^S$ corresponding to $j\in S_1$ and set them to $\epsilon$. 
\hfill \QED
\subsection{Proof of Proposition \ref{thm_submodularinfluence}}

% \begin{theorem} \label{thm_submodularinfluence}
% The influence function $f_v(\cdot)$ is almost surely submodular in % dense networks.
% \end{theorem}

To show submodularity of the set function $f_v(S)$, we need to show that for all $S \subseteq S' \subset S_0$ and $u \in S_0\setminus S', u\neq v$, we have 
\begin{align}
    f_v(S\cup \{u\})-f_v(S) \geq f_v(S' \cup \{u\})-f_v(S')  \label{eq_subinf}
\end{align}
Consider the first term in \eqref{eq_subinf}. 
\begin{align}
    f_v(S\cup \{u\})&=g^{v}\mathbb{1}(1_n-diag(\textbf{c})\Mat{G}\lceil a^{S\cup \{u\}} \rceil) \nonumber \\
    &= g^{v}\mathbb{1}\bigg(1_n - \big(diag(\textbf{c})\Mat{G}\lceil a^{S} \rceil - diag(\textbf{c})(g^{u})^T\big)\bigg)
    \label{eq_infbreakdown}
\end{align}
The additional term in \eqref{eq_infbreakdown} arrives from the fact that if $u \in S_0$ is additionally part of the zero-set, then all $j\in n(u)$ have their neighbor influence reduced by $c_j$. The non-neighbors of $u$ in set $S_1$ are unaffected. Now let
\begin{equation}
\begin{aligned}
  A&=\mathbb{1}(\bbone_n - diag(\textbf{c})\Mat{G}\lceil a^{S} \rceil + diag(\textbf{c})(g^{u})^T) \\
  B&=\mathbb{1}(\bbone_n-diag(\textbf{c})\Mat{G}\lceil a^{S} \rceil)  \\
  Z&=\mathbb{1}(\bbone_n-diag(\textbf{c})\Mat{G}\lceil a^{S'} \rceil)  \\
  D&=\mathbb{1}(\bbone_n - diag(\textbf{c})\Mat{G}\lceil a^{S'} \rceil + diag(\textbf{c})(g^{u})^T)
\end{aligned}
\end{equation}

% $A=\mathbb{1}(1_N - diag(\textbf{c})\Mat{G}\lceil a^{S} \rceil + diag(\textbf{c})(g^{u})^T)$, $B=\mathbb{1}(1_N-diag(\textbf{c})\Mat{G}\lceil a^{S} \rceil)$, $Z=\mathbb{1}(1_N-diag(\textbf{c})\Mat{G}\lceil a^{S'} \rceil)$, $D=\theta(1_N - diag(\textbf{c})\Mat{G}\lceil a^{S'} \rceil + diag(\textbf{c})(g^{u})^T)$. 
% \\

\textbf{\textit{Case 1}}: Consider $w\in n(v)$, 
% \blue{CE: you cannot use $n$ here. $n$ is the number of players and the neighbor function...}, 
$w\not\in n(u)$. Then
\begin{align}
    A_w-B_w+Z_w=\mathbb{1}(1-c_w g^{w} \lceil a^{S} \rceil) - \mathbb{1}(1-c_w g^{w} \lceil a^{S} \rceil) \nonumber \\
    +\mathbb{1}(1-c_w g^{w} \lceil a^{S'} \rceil) = D_w
\end{align}
which implies
\begin{align}
    \sum_{w\in n(v)} A_w -B_w + Z_w = \sum_{w \in n(v)} D_w \nonumber
\end{align}
or,
\begin{align}
    g^{v}(A-B+Z)&=g^{v}(D) \nonumber \\
    \implies f_v(S\cup \{u\})-f_v(S) &= f_v(S' \cup \{u\})-f_v(S') 
\end{align}
\textbf{\textit{Case 2}}: Take $w\in n(v)$ $\cap$ $n(u)$.
% \blue{CE: Same here. I'd choose another letter to denote the agent...}. 
Consider now $\lceil a^S \rceil$ and $\lceil a^{S'} \rceil$. For all $i \in S_1$, $\lceil a^S_i \rceil=\lceil a^{S'}_i \rceil$. Since $S\subseteq S'$, $\forall i \in S_0$, $\lceil a^S_i \rceil \geq \lceil a^{S'}_i \rceil$. We shall hereby use the following ordering for simplicity, 
\begin{align}
    \lceil a^S \rceil \geq \lceil a^{S'} \rceil \label{eq_vecOrd}
\end{align}
Now it is easy to show that the following inequalities hold
\begin{align}
    B_w \leq Z_w \leq D_w \nonumber \\
    B_w \leq A_w \leq D_w \label{eq_infconstr}
\end{align}
using $c_w g^{w} \lceil a^S \rceil \geq c_w g^{w} \lceil a^{S'} \rceil$ which follows from \eqref{eq_vecOrd}. For condition \eqref{eq_subinf} to be violated, $\exists$ $w \in n(v)$ such that
\begin{align}
    D_w > A_w -B_w + Z_w \label{eq_subinfViol}
\end{align}
holds. Next we analyze the chance that the above inequality is realized. From the inequalities in \eqref{eq_infconstr}, we see that \eqref{eq_subinfViol} holds only when $D_n=1$ and $A_n=B_n=Z_n=0$ (see Table (\ref{tab:tableABCD})). For that to happen, we need 
\begin{align}
  -c_w &< 1-c_w g^{w} \lceil a^{S'} \rceil \leq 0 \label{eq_ineq1},
\end{align}
and
\begin{align}
   1-c_w g^{w} \lceil a^{S} \rceil &\leq -c_w. \label{eq_ineq2}
\end{align}
% and 
% \begin{align}
%    
% \end{align}
We rewrite the inequalities in \eqref{eq_ineq1} and \eqref{eq_ineq2} respectively as
\begin{align}
    \frac{1}{g^{w}\lceil a^{S'} \rceil}&\leq c_w<\frac{1}{g^{w}\lceil a^{S'} \rceil - 1}
\end{align}
\begin{align}
    \frac{1}{g^{w}\lceil a^{S} \rceil-1}&\leq c_w.  \label{eq_ineq2A}
\end{align}
Since $\lceil a^S \rceil \geq \lceil a^{S'} \rceil$, we have
\begin{align}
    \frac{1}{g^{w}\lceil a^S \rceil} \leq \frac{1}{g^{w}\lceil a^{S'} \rceil}
\end{align}
 and 
 \begin{align}
    \frac{1}{g^{w}\lceil a^S \rceil-1} \leq \frac{1}{g^{w}\lceil a^{S'} \rceil-1} 
 \end{align}
  respectively. Now two cases can arise. For \eqref{eq_subinf} to fail, that is for \eqref{eq_subinfViol} to be satisfied, we must have $c_w$ falling in the range 
\begin{align}\label{eq_range0}
    \frac{1}{g^{w}\lceil a^S \rceil-1} \leq c_w\leq \frac{1}{g^{w}\lceil a^{S'} \rceil-1}, \\ \nonumber
    \textrm{whenever} \quad \frac{1}{g^{w}\lceil a^S \rceil-1} \geq \frac{1}{g^{w}\lceil a^{S'} \rceil}  
\end{align}
is satisfied, or
\begin{align}
     \frac{1}{g^{w}\lceil a^{S'} \rceil} \leq c_w \leq \frac{1}{g^{w}\lceil a^{S'} \rceil-1}, \textrm{otherwise.} \label{eq_range}
\end{align}
% \blue{CE: In (53), should it not be $S'$ in the first fraction?} Nope
Out of the two ranges in \eqref{eq_range0} and \eqref{eq_range}, the latter is bigger.
We denote the probability of $c_w$ falling in the range in \eqref{eq_range} by $p_{c_w}$. We have the following bound 
\begin{align}
    p_{c_w} \leq \big(\frac{1}{g^{w}\lceil a^{S'} \rceil-1} - \frac{1}{g^{w}\lceil a^{S'}\rceil}\big)\frac{d_w}{d_w-1}. \label{eq_cnrange}
\end{align}
Here $d_w := |n(w)|$ is the degree of agent $w$.
% \blue{CE: right after the equation you need to define $d_n$. Also it is not clear how you got the second fraction...}

This is because $c_w$ is uniformly distributed between $[\frac{1}{d_n},1)$ (since we have a well-behaved instance almost surely whenever the network is dense). See that for a uniformly distributed random variable in the range $[a,b]$, the probability of it lying in the range $[a',b'], a<a'<b'<b$ is given by $(b'-a')/(b-a)$.

\begin{table}[t]
\centering
\begin{tabular}{c|c|c|c|c|c|c}   \toprule
$A_w -B_w +Z_w$ & $1$ & $2$ &$2$ & $1$ &$0$ &$0$\\ %& $C(3,0)$ &$C(3,1)$ & $C(3,2)$ &$C(3,3)$\\
\hline
$D_w$ & 0 & 0& 1& 1& 0& 1\\%& 3 & 4 & 4& 6\\

\bottomrule
\end{tabular}
\caption{ Possible combination of values for $A_w -B_w +Z_w$ and $D_w$}
\label{tab:tableABCD}
\vspace{-12pt}
\end{table}
Further simplifying \eqref{eq_cnrange} and noting that $g^{w}\lceil a^{S'} \rceil \leq d_w$, we see that 
\begin{align}
    p_{c_w} \leq O(\frac{1}{d_{w}^2}).
\end{align}

Since, for dense networks $\frac{1}{d_{w}^2} \rightarrow 0$, we have the desired result.
\hfill \QED
\begin{corollary} \label{cor_monotontinfluence}
The one-step influence function $f_v(\cdot)$ is monotone increasing.
\end{corollary}
\begin{proof}
 For $S\subset S^{'} \subset S_0$ we have $\lceil a^S \rceil \geq \lceil a^{S'} \rceil$ which implies that
\begin{align}
 \bbone_n - diag(\textbf{c})\Mat{G}\lceil a^S \rceil \leq \bbone_n-diag(\textbf{c})\Mat{G}\lceil a^{S'} \rceil    
\end{align}
 which further implies the result.
\end{proof}

\subsection{Proof of Lemma \ref{lem_process}} \label{sec_lemma_process}
Let $E(A_t)$ be the set of edges deactivated at time $t$ when we initiate the dynamics in \eqref{eq_updates} with control set $A$. Since we have a well-behaved instance almost surely, at any time $t$, we have some agents in $S_0$ playing $0$, some agents in $S_1$ playing $1$, the rest undecided. This implies that all edges deactivated at time $t$ are incident with the zero-set at time $t$. Recall that we have finite time convergence of the dynamics. As a consequence of the aforementioned coupling, we have $Z_{\infty} \subseteq A_{\infty} \cap B_{\infty}$ whereby the set of deactivated incident edges $E(Z_{\infty}) \subseteq E(A_{\infty}) \cap E(B_{\infty})$. Moreover $D_{\infty} \subseteq A_{\infty} \cup B_{\infty}$ implies $E(D_{\infty}) \subseteq E(A_{\infty}) \cup E(B_{\infty})$. Recall we can represent the objective in MAC \eqref{eq_MAC} as $f(X)=|E(X_{\infty})|$. We can now write
\begin{align}
    f(A)+ f(B) &= |E(A_{\infty})|+|E(B_{\infty})| \nonumber \\
    &= |E(A_{\infty})\cap E(B_{\infty})| + |E(A_{\infty})\cup E(B_{\infty})| \nonumber \\
    &\geq |E(Z_{\infty})|+|E(D_{\infty})|\nonumber \\
    &=f(Z)+f(D)\nonumber \\
    &=f(A\cap B) + f(A\cup B)
\end{align}
The second equality is due to the fact that the cardinality function is both submodular and supermodular, which means that for any two sets $I_1,I_2$ subset of a ground set $I_0$ we have 
\begin{align}
    |I_1|+|I_2|=|I_1 \cup I_2|+|I_1 \cap I_2|
\end{align}
The first inequality is a result of monotonicity of $|\cdot|$.
Taking expectation on either side, we get \eqref{subexp}.
\hfill \QED

\bibliographystyle{IEEEtran}
\bibliography{IP}
\end{document}